\title{\vspace{-1cm} Avoiding small subgraphs in Achlioptas processes}

\author{
Michael Krivelevich \thanks{School of Mathematical Sciences,
Raymond and Beverly Sackler Faculty of Exact Sciences, Tel Aviv
University, Tel Aviv 69978, Israel. E-mail:
krivelev@post.tau.ac.il.
Research supported in part by USA-Israel BSF Grants 2002-133 and
2006-322, and
by grant 526/05 from the Israel Science Foundation.}
\and 
Po-Shen Loh \thanks{Department of Mathematics,
Princeton University, Princeton, NJ 08544. E-mail:
ploh@math.princeton.edu.
Research supported in part by a Fannie and John Hertz Foundation Fellowship, an NSF Graduate
Research Fellowship, and a Princeton Centennial Fellowship.}
\and 
Benny Sudakov \thanks{Department of Mathematics, UCLA, Los Angeles, CA 90095. E-mail:
bsudakov@math.ucla.edu.
Research supported in part by NSF CAREER award DMS-0812005, and a USA-Israeli BSF grant.
}
}

\documentclass[11pt]{article}
\usepackage{amsmath, amssymb, amsthm, epsfig, graphicx, enumerate}

\newtheorem{theorem}{Theorem}[section]
\newtheorem*{theorem*}{Theorem}

\newtheorem{lemma}[theorem]{Lemma}
\newtheorem{inequality}[theorem]{Inequality}
\newtheorem*{lemma*}{Lemma}

\newtheorem{corollary}[theorem]{Corollary}

\newtheorem{definition}[theorem]{Definition}

\newcommand{\ep}{\epsilon}
\newcommand{\pr}[1]{\mathbb{P}\left[#1\right]}

\newcommand{\bin}{\text{Bin}}

\newcommand{\whp}{\textbf{whp}}
\newcommand{\wep}{\textbf{wep}}
\newcommand{\hme}[2]{#1 \setminus #2e}

\newcommand{\shortversion}[2]{#2}  

\date{}
\begin{document}
\maketitle

\begin{abstract}
  For a fixed integer $r$, consider the following random process.  At
  each round, one is presented with $r$ random edges from the edge set
  of the complete graph on $n$ vertices, and is asked to choose one of
  them.  The selected edges are collected into a graph, which thus grows
  at the rate of one edge per round.  This is a natural generalization
  of what is known in the literature as an \emph{Achlioptas process}\/
  (the original version has $r=2$), which has been studied by many
  researchers, mainly in the context of delaying or accelerating the
  appearance of the giant component.

  In this paper, we investigate the \emph{small subgraph}\/ problem
  for Achlioptas processes.  That is, given a fixed graph $H$, we
  study whether there is an online algorithm that substantially delays or
  accelerates a typical appearance of $H$, compared to its threshold of
  appearance in the random graph $G(n, M)$.  It is easy to see that
  one cannot accelerate the appearance of any fixed graph by more than
  the constant factor $r$, so we concentrate on the task of avoiding
  $H$.  We determine thresholds for the avoidance of all cycles $C_t$,
  cliques $K_t$, and complete bipartite graphs $K_{t,t}$, in every
  Achlioptas process with parameter $r \geq 2$.
\end{abstract}

\section{Introduction}

The standard Erd\H{o}s-R\'enyi random graph model $G(n, M)$ can be
described as follows.
Start with the empty graph on $n$ vertices, and perform $M$ rounds, adding one random edge
to the graph at each round.  For any monotone increasing graph property (such as
containment of $K_4$ as a subgraph, say), it is natural to ask whether there is some value
of $M$ at which the probability of $G(n, M)$ satisfying the property changes rapidly from
nearly 0 to nearly 1. More precisely, a function $M^*(n)$ is said to be a threshold for a
property $\cal P$ if for any $M(n) \ll M^*(n)$, the random graph $G(n, M)$ does not
satisfy $\cal P$ \whp, but for any $M(n) \gg M^*(n)$, the random graph $G(n, M)$ satisfies
$\cal P$ \whp.  Here, \whp\ stands for \emph{with high probability}, that is, with
probability tending to 1 as $n \rightarrow \infty$, and $f(n) \ll g(n)$ means that $f(n)/g(n) \rightarrow 0$ as $n \rightarrow \infty$. A classical result of Bollob\'as and
Thomason \cite{BT} implies that every monotone graph property has a threshold, and much
work has been done to determine thresholds for various properties.

Recently, there was much interest in the following natural variant of
the classical model.  We still begin with the empty graph and
perform a series of rounds, but at each round, one is now presented
with two independent and uniformly random edges, and is asked to
choose one of them to add to the graph.  This is known in the
literature as an \emph{Achlioptas process}, after Dimitris Achlioptas,
who asked the question of whether there was an online algorithm which could,
with high probability, substantially delay the appearance of the giant
component (a connected component with $\Omega(n)$ vertices).

The
trivial algorithm, which arbitrarily chooses the first edge in each offered pair, 
essentially produces
the random graph $G(n, M)$ after $M$ rounds, so $G(n, M)$ serves as the
benchmark against which comparisons are made.
A classical result of Erd\H{o}s and R\'enyi \cite{ER} states that if
$M = cn$ for any absolute constant $c > 1/2$, then the random graph
$G(n, M)$ contains a giant component \whp.  For the Achlioptas process,
Bohman and Frieze \cite{BF} presented an algorithm which could run for $0.535n$ rounds, 
while keeping 
the size of the largest component only poly-logarithmic in $n$ \whp.  
Since then, much work has been done \cite{BBFP1, BBFP2, BFW, BKi, BKr, FGS, SW}.  
The current best result
for this problem is due to Spencer and Wormald \cite{SW}, who
exhibit an algorithm that can run for $0.829n$ rounds while
keeping all component sizes bounded by $O(\log n)$ \whp.  In the opposite
direction, Bohman, Frieze, and Wormald \cite{BFW} have shown that no
algorithm can succeed \whp\ past $0.964n$ rounds.  Several variants
have also been studied, such as the offline version, a two-player
version, and the question of \emph{embracing}\/ (accelerating the
appearance of) the giant component.

While the main focus of the research mentioned above was the giant
component, it is natural to study other graph properties in the
context of Achlioptas processes.  In this paper, we study the problem
which in the literature is referred to as the \emph{small subgraph}\/
problem.  This was one of the main problems studied in the seminal
paper of Erd\H{o}s and R\'enyi \cite{ER} from 1960, which was the
starting point of the theory of random graphs.  The original problem,
stated for the random graph model $G(n, M)$, was as follows: given a
fixed graph $H$ (a triangle or $K_4$, say), find the smallest value of
$M$ such that the random graph $G(n, M)$ contains $H$ as a (not
necessarily induced) subgraph \whp.  The subgraph is called ``small''
because its size is fixed while $n$ tends to infinity.

It turns out that in this problem, the relevant parameter is the
maximum edge density $m(G) = \max\{e(H)/v(H) : \text{$H$ is a subgraph
  of $G$} \}$.  In their original paper, Erd\H{o}s and R\'enyi found
thresholds for all balanced graphs, which are the graphs whose edge
density $e(H)/v(H)$ equals the maximum edge density $m(H)$.  It
was not until 20 years later that Bollob\'as \cite{B} solved the problem for
all graphs, proving that for any $H$ with $m(H) \geq 1$, the threshold
for $H$ appearing in $G(n, M)$ is $M^* = n^{2 - \frac{1}{m(H)}}$.  For
further reading about the small subgraph problem in $G(n, M)$, we
direct the interested reader to the monographs by Bollob\'as
\cite{B-book} and by Janson, \L uczak, and Ruci\'nski \cite{JLR}, each
of which contains an entire section discussing this problem.

In this paper, we consider the small subgraph problem in the context
of Achlioptas processes, and investigate whether one can substantially
affect thresholds by introducing this power of choice.  Actually, we
study a natural generalization of the process, which we call an
\emph{Achlioptas process with parameter $r$}.  In this process, $r$
edges of $K_n$ are presented at each round, and one of them is
selected.  We will always consider $r$ to be fixed as $n$ tends to
infinity (note that $r=2$ corresponds to the original Achlioptas
process).

Let us now state our model precisely.  At the $i$-th round, one is
presented with $r$ independent random edges, each distributed
uniformly over all ${n \choose 2} - (i-1)$ remaining edges that have
not yet been chosen for the graph.  Note that this eliminates the
possibility of choosing the same edge twice, so our final graph is
simple.  However, we do allow the possibility that edges may be
offered more than once, which simplifies our arguments.  One may
consider models in which all sampling is with replacement (which may
create multigraphs), or in which every edge is offered at most once,
but our results in this paper will still carry over because we always
run the process for $o(n^2)$ rounds.

Note that the graph after the $k$-th round of the Achlioptas process with
parameter $r$ is a subgraph of the random graph with $rk$ edges.  So,
the question of accelerating the appearance of a fixed graph is
immediately resolved in the negative.  Clearly, the threshold cannot move
forward by more than a (constant) factor of $r$.

So, in this paper we concentrate on the avoidance problem.  We may
pose it as a single player game in which the player loses when he
creates a (not necessarily induced) subgraph isomorphic to a certain
fixed graph $H$.  The player's objective is to postpone losing for as
long as possible.  We say that a function $m^*(n)$ is a
\emph{threshold for avoiding $H$}\/ if: \textbf{(i)} given any
function $m(n) \ll m^*(n)$, there exists an online strategy by which the
player survives through $m$ rounds \whp, 
and \textbf{(ii)} given any function $m(n)
\gg m^*(n)$, the player loses by the end of $m$ rounds \whp, regardless
of the choice of such a strategy.

Note, however, that it is not obvious that thresholds necessarily
exist.  Furthermore, unlike the situation in the small subgraph
problem, there are no simple first-moment calculations that suggest
what the thresholds should be.  As it turned out, a substantial part
of the difficulty in obtaining our results was in conjecturing the
correct thresholds.  We were able to solve the problem for all cycles
$C_t$, cliques $K_t$, and complete bipartite graphs $K_{t,t}$.  Let us
now state our main result:

\begin{theorem} \mbox{}
  \label{thm:main}

  \begin{description}

  \item[(i)] For $t \geq 3$, the threshold for avoiding $C_t$ in the
    Achlioptas process with parameter $r \geq 2$ is
    $n^{2-\frac{(t-2)r+2}{(t-1)r+1}}$.

  \item[(ii)] For $t \geq 4$, the threshold for avoiding $K_t$ in the
    Achlioptas process with parameter $r \geq 2$ is $n^{2-\theta}$,
    where $\theta$ is defined as follows:
    \begin{displaymath}
      s = \lfloor \log_r [(r-1)t + 1] \rfloor, \quad\quad\quad \theta = \frac{r^s (t-2) + 2}{r^s\left[{t \choose 2} - s\right] + \frac{r^s - 1}{r-1}}.
    \end{displaymath}

  \item[(iii)] For $t \geq 3$, the threshold for avoiding $K_{t,t}$ in
    the Achlioptas process with parameter $r \geq 2$ is $n^{2-\theta}$, where
    $\theta$ is defined as follows:
    \begin{displaymath}
      s = \lfloor \log_r [(r-1)t + 1] \rfloor, \quad\quad\quad \theta = \frac{r^s (2t-2) + 2}{r^s (t^2 - s) + \frac{r^s - 1}{r-1}}.
    \end{displaymath}

  \end{description}
  
\end{theorem}


\noindent \textbf{Remark.}\, In all of these cases, we provide
\emph{deterministic}\/ online algorithms that achieve the thresholds
\whp, but show that even randomized algorithms cannot survive
beyond the thresholds.

\vspace{4mm}

The rest of this paper is organized as follows.  In the next section, we present
some tools from extremal combinatorics and the theory of random
graphs, which we will use in our proofs.  Then, we present the proof
of our theorem, which is divided into several sections.  We begin in
Section \ref{sec:warm-up} with the case of avoiding $K_4$ when $r=2$,
which turns out to be the first
nontrivial case.  We treat this case in detail, because our argument
there is the prototype for the general argument that we later use to
prove thresholds for $K_t$, $K_{t,t}$, and $C_t$.

We extend the argument to almost all other $K_t$ and $r$ in Section
\ref{sec:avoid-kt}.
The proof requires many
inequalities whose somewhat tedious verifications would interfere with
the exposition,
\shortversion{%
  so their precise statements are recorded in the
  appendix.\footnote{The proofs of these inequalities are rather technical
    and not so interesting, so they only appear in the unabridged
    version of this paper which is on the arXiv at
    \texttt{http://arxiv.org/abs/0708.0443}.}
}%
{%
  so we postpone their proofs to the appendix.
}%
This also makes it easier to distill the abstract argument,
which we present in Section \ref{sec:abstract-argument}.  Next, we
apply the abstraction to prove thresholds for avoiding $C_t$ in
Section \ref{sec:avoid-cycles} and $K_{t,t}$ in Section \ref{sec:avoid-ktt}.
We treat the last remaining case of
avoiding $K_4$ in the Achlioptas process with parameter 3 in Section
\ref{sec:avoid,t=4,r=3}.  The final section contains some concluding
remarks and open problems.

\section{Preliminaries}

\subsection{Notation and terminology}
\label{sec:notation}

Throughout our paper, we will omit floor and ceiling signs whenever
they are not essential, to improve clarity of presentation.  The
following (standard) asymptotic notation will be utilized extensively.
For two functions $f(n)$ and $g(n)$, we write $f(n) = o(g(n))$ or
$g(n) = \omega(f(n))$ if $\lim_{n \rightarrow \infty} f(n)/g(n) = 0$,
and $f(n) = O(g(n))$ or $g(n) = \Omega(f(n))$ if there exists a
constant $M$ such that $|f(n)| \leq M|g(n)|$ for all sufficiently
large $n$.  We also write $f(n) = \Theta(g(n))$ if both $f(n) =
O(g(n))$ and $f(n) = \Omega(g(n))$ are satisfied.

Let us introduce the following abbreviations for some phrases that we
will use many times in our proof.  As mentioned in the introduction,
\whp\ will stand for \emph{with high probability}, i.e., with
probability $1 - o(1)$.  It is also convenient for us to introduce the
abbreviation \wep, which stands for \emph{with exponential
  probability}, i.e., with probability $1 - o\big(e^{-n^c}\big)$ for
some $c>0$.  We will say that a function $f$ is a \emph{positive power of
  $n$}\/ if  $f = \Omega(n^c)$ for some $c>0$.  Analogously, we will
say that a function $f$ is a \emph{negative power of $n$}\/ if 
$f = O(n^{-c})$ for some $c>0$.

Next, let us discuss the graph-specific terms that we will use.  We
often need to consider the graphs at intermediate stages of the
Achlioptas process, so $G_i$ will always denote the graph after the
$i$-th round.  Our main interest in $G_i$ will be to count
\emph{copies}\/ of subgraphs.  Here, we define a copy of a graph $H$
in another graph $G$ to be an injective map from $V(H)$ to $V(G)$ that
preserves the edges of $H$.  Note that copies are
not necessarily induced subgraphs, and are labeled, i.e.,  we do not take automorphisms
into account when computing the number of copies of $H$ in a graph.

The player's objective in the Achlioptas process is to avoid creating
a copy of a certain fixed graph $H$, but our analysis needs to
consider subgraphs of $H$ as well.  It is therefore convenient to introduce
the notation $\hme{H}{k}$ to represent any graph which can be obtained
by deleting any $k$ edges from $H$.  (When $k = 1$, we will simply
write $\hme{H}{}$.)  This enables us to concisely refer to all graphs
of the form $\hme{H}{k}$ in the aggregate.  For example, the phrase
``the number of copies of $\hme{H}{k}$'' should be understood to be the
total number of copies of all graphs of the form $\hme{H}{k}$.

We keep track of the numbers of copies of these subgraphs by studying how counts are
affected by the addition of an edge at a pair of vertices.  This motivates the following
definition.  Let $G$ and $H$ be graphs, let $k$ be an integer, and let $a, b$ be a pair of
distinct vertices of $G$.  Let $G^+$ be the graph obtained from $G$ by adding the edge
between $a$ and $b$ if it is not yet present, and let $G^-$ be the graph obtained by
deleting that edge if it was present.  Note that $G$ is equal to either $G^+$ or $G^-$.
Then, we say that the pair $\{a, b\}$ \emph{completes $t$ copies of $\hme{H}{k}$}\/ if $t$
is the difference between the number of copies of $\hme{H}{k}$ in $G^+$ and the number in
$G^-$.

Sometimes, we need to be specific about which graphs of the form
$\hme{H}{(k+1)}$ are completed into graphs of the form $\hme{H}{k}$.
Let $H_1$ and $H_2$ be graphs on the same vertex set $U$, with $E(H_1)
\subset E(H_2)$, but differing only in exactly one edge.  Let $\{u,
v\} \subset U$ be the endpoints of that edge.  Let $G$ be another
graph, and let $a, b$ be a pair of distinct vertices of $G$.  Then, we
say that the pair $\{a, b\}$ \emph{extends $t$ copies of $H_1$ into
  $H_2$}\/ if $t$ is the number of injective graph homomorphisms $\phi : H_1
\rightarrow G$ that map $\{u, v\}$ to $\{a, b\}$.
Note that this definition is insensitive to the presence of an edge
between $a$ and $b$.

\subsection{Extremal combinatorics}

In this section, we present two extremal results, which are used
in the proofs of the upper bounds in our thresholds (i.e., that no
strategy can survive for too many rounds).  The
following lower bound on the number of paths in a graph was
obtained in \cite{ES} using a matrix inequality of Blackley and Roy.

\begin{lemma}
  \label{lem:count-paths}
  Every graph with $n$ vertices and average degree $d$ contains at
  least $(1+o(1)) n d^{t-1}$ copies of the $t$-vertex path $P_t$.
  Here, we consider $t$ to be fixed, while $d$ and $n$ tend to
  infinity.
\end{lemma}

Next, we record the following well-known extremal result, which lower
bounds the number of copies of the complete bipartite graph $K_{s,t}$
that can appear in any graph with a fixed number of edges.  The
classical proof (via two applications of convexity) is based on the
ideas used by K\"ov\'ari, S\'os, and Tur\'an \cite{KST} to bound
the Tur\'an number $\text{ex}(n, K_{s,t})$.

\begin{lemma}
  \label{lem:extremal-count-ktt}
  For fixed positive integers $s \leq t$, and any function $p \gg
  n^{-1/s}$, every graph with $n$ vertices and ${n \choose 2} p$ edges
  contains at least $(1+o(1)) n^{s+t} p^{st}$ copies of the complete
  bipartite graph $K_{s, t}$.
\end{lemma}

\subsection{Random graphs}
\label{sec:random-graphs}

We begin by recalling the Chernoff bound for exponential
concentration of a binomial random variable.  We use the formulation
from \cite{AS}.

\begin{theorem} 
  For any $\ep > 0$, there exists $c_\ep > 0$ such that
    the following holds.  Let $X$ be any binomial random variable, and
    let $\mu$ be its expectation.  Then $\pr{| X - \mu |  > \ep \mu} < 2 e^{-c_\ep \mu}$.
\end{theorem}

Using the Chernoff bound and a standard coupling argument, we prove a
result that allows us to relate $G_m$ (the graph after the $m$-th
round of the Achlioptas process) to the more familiar random graph
$G(n, p)$.

\begin{lemma}
  \label{lem:coupling}
  Suppose that
  $n \ll m \ll n^2$.  Then we may couple the Achlioptas process with $G(n, p =
  4rm/n^2)$ in such a way that \wep, $G_m$ is a subgraph of $G(n, p)$.
\end{lemma}

\noindent {\bf Proof.}\, In
the Achlioptas process, $r$ random edges are presented at each round,
independently and uniformly distributed over all potential edges that
have not yet been picked for the graph.  So, we may couple the first
$m$ rounds of the process with the edge-uniform random graph $G(n,
rm)$ in such a way that if we consider the graph $G_m^+$ obtained by
taking every edge that was offered (instead of choosing only one per
round), $G_m^+$ is always a subgraph of $G(n, rm)$.  Yet $G_m$ is
always a subgraph of $G_m^+$, so it remains to relate $G(n, rm)$ with
$G(n, p=4rm/n^2)$.  This final part is standard and proceeds via
coupling with the random graph process; under this coupling, $G(n, rm)
\subset G(n, p)$ as long as $\bin\big[{n \choose 2}, p\big] \geq rm$,
and the Chernoff bound shows that this event occurs \wep.  \hfill
$\Box$

\vspace{4mm}

Our analysis revolves around counting copies of fixed subgraphs in
$G_m$.  The previous lemma allows us to apply results from the theory
of $G(n, p)$ to assist us in this pursuit.  We now record several such
theorems, translated in terms of $G_m$.  The following definition is
crucial for counting subgraphs in $G(n, p)$.

\begin{definition}
  A graph $H$ is \textbf{balanced} if for any subgraph
  $H'\subseteq H$, $\frac{e(H')}{v(H')} \leq
  \frac{e(H)}{v(H)}$.
\end{definition}

\begin{theorem}
  \label{thm:balanced-graph}
  Let $H$ be a fixed balanced graph with $v$ vertices and $e$ edges.
  Suppose that
  $n \ll m \ll n^2$, and let $p = 2m/n^2$.  Also suppose that $n^v
  p^e$ is a positive power of $n$.  Then the number of copies of $H$
  in $G_m$ is $O(n^v p^e)$ \wep.
\end{theorem}

\noindent {\bf Proof.}\, By Lemma \ref{lem:coupling}, it suffices to
count copies of $H$ in $G(n, 2rp)$.  The expected number of copies is
$(1+o(1)) n^v (2rp)^e = \Theta(n^v p^e)$, which is a positive power of $n$ by
assumption.  This allows us to apply Corollary 6.3 of \cite{V}, which
uses Kim-Vu polynomial concentration to prove the following result:
for any balanced graph $H$ such that the expected number of copies of
$H$ in the random graph is $\mu \gg \log n$, the probability that the
actual number of copies exceeds $2\mu$ is $e^{-\Omega(\mu)}$.  In
our case, $\mu$ is a positive power of $n$, so this implies that \wep,
the number of copies is $O(n^v p^e)$, as desired.  \hfill $\Box$

\vspace{4mm}

The previous result provides a very precise count of the number of
copies of a fixed graph in the random graph $G(n, p)$.  However, the
point of the Achlioptas process was to deviate from $G(n, p)$ by
introducing the power of choice.  So, our analysis will have to take
the potential of choice into account.  We keep track of the numbers of
copies of subgraphs by studying how counts are affected by the
addition of an edge at a pair of vertices; this motivated the notions
of a pair \emph{completing $t$ copies of $\hme{H}{k}$}\/ and of the
pair \emph{extending $t$ copies of $H_1$ into $H_2$}, which we defined
at the end of Section \ref{sec:notation}.  

This is essentially the problem of counting extensions, which has also
been well-studied in $G(n, p)$.  We refer the interested reader to
Chapter 10 of \cite{AS}.  As in the case of counting subgraphs in
$G(n, p)$, a suitable definition of balanced-ness is required to count
extensions.

\begin{definition} \mbox{}
  \begin{description}
  \item[(i)] Let $H_1$ and $H_2$ be graphs on the same vertex set $U$,
    with $E(H_1) \subset E(H_2)$, but differing only on the edge
    joining the vertices $u, v \in U$.  We say that the pair $(H_1,
    H_2)$ is a \textbf{balanced extension pair} if for every
    proper subset $U' \subset U$ that still contains $\{u, v\}$, the
    induced subgraph $H' = H_1[U']$ has the property that
    $\frac{e(H')}{v(H')-2} \leq \frac{e(H_1)}{v(H_1)-2}$.

  \item[(ii)] $\hme{H}{k}$ has the \textbf{balanced extension property} if
    every pair $(H_1, H_2)$ with $V(H_1) = V(H_2) = V(H)$, $E(H_1) \subset E(H_2) \subset
    E(H)$, $e(H_1) = e(H) - k$, and $e(H_2) = e(H) - k + 1$, is a balanced extension
    pair.
  \end{description}
\end{definition}

\begin{theorem}
  \label{thm:balanced-extension}
  Suppose that
  $n \ll m \ll n^2$, and let $p = 2m/n^2$.  Let $(H_1, H_2)$ be a
  balanced extension pair, and let $v$ and $e$ be the numbers of
  vertices and edges in $H_1$, respectively.  Finally, let $j$ be an
  arbitrary integer constant.
  \begin{description}
  \item[(i)] Suppose that $n^{v-2} p^e$ is a positive power of $n$.
    Then \wep, every pair of distinct vertices $\{a, b\}$ of $G_{jm}$
    extends $O(n^{v-2} p^e)$ copies of $H_1$ into $H_2$.
  \item[(ii)] Suppose that $n^{v-2} p^e$ is a negative power of $n$.  Then, for any
    constant $\gamma > 0$, there exists a constant $C$ such that with probability $1 -
    o(n^{-\gamma})$, every pair of distinct vertices $\{a, b\}$ of $G_{jm}$ extends at most $C$
    copies of $H_1$ into $H_2$.
  \end{description}
\end{theorem}

\noindent {\bf Proof.}\, By Lemma \ref{lem:coupling}, it suffices to
consider $G(n, 2rjp)$ instead of $G_{jm}$ in both parts of the
theorem.  For part (i), the expected number of extensions at a pair in
$G(n, 2rjp)$ is $(1+o(1)) n^{v-2} (2rjp)^e = \Theta(n^{v-2} p^e)$,
which is a positive power of $n$ by assumption.  This allows us to
apply Corollary 6.7 of \cite{V}, which uses Kim-Vu polynomial
concentration to prove the following result: for any balanced
extension pair $(H_1, H_2)$ such that the expected number $\mu$ of
copies of $H_1$ that a fixed edge extends into $H_2$ in the random
graph is a positive power of $n$, the probability that the actual
number of extensions exceeds $2\mu$ is $e^{-\Omega(\mu)}$.  In
our case, $\mu$ is a positive power of $n$, so even after taking a
union bound over all $O(n^2)$ pairs of vertices, this implies that
\wep, every pair of vertices extends $O(n^{v-2} p^e)$ copies of
$H_1$ into $H_2$.  This establishes (i).

For part (ii), let us bound the probability that $\{a, b\}$ extends
$C$ copies of $H_1$ into $H_2$.  Recall that $H_1$ and $H_2$ shared
the same vertex set $U$, and differed only on the edge joining $u, v
\in U$.  Consider any graph $F$ which is formed by the superposition
of $C$ distinct copies of $H_1$, all with $\{u, v\}$ mapping to the
same pair of vertices $\{u', v'\} \in V(F)$.  Let $v' = v(F)$ and $e'
= e(F)$.  

The probability that $\{a, b\}$ has an extension to $F$ (an injective
map from $V(F)$ sending $\{u', v'\} \mapsto \{a, b\}$) in $G(n, 2rjp)$
is at most $n^{v'-2} (2rjp)^{e'} = O((n p^{e'/(v'-2)})^{v'-2})$.  An easy and
standard induction, using the fact that $(H_1, H_2)$ is a balanced
extension pair, implies that $\frac{e'}{v'-2} \geq \frac{e}{v-2}$.
Hence this probability is at most $O((n p^{e/(v-2)})^{v'-2}) = O((n^{v-2}
p^e)^\frac{v'-2}{v-2})$.

We assumed that $n^{v-2} p^e$ was a negative power of $n$.  Also,
since the $C$ copies of $H_1$ in $F$ are distinct, one can trivially
bound $C \leq (v'-2)^{v-2} \Rightarrow v'-2 \geq C^\frac{1}{v-2}$.
So, for a sufficiently large constant $C$, the probability that $\{a,
b\}$ has an extension to $F$ is $o(n^{-\gamma-2})$.  Taking a union
bound over all $O(n^2)$ pairs of vertices, we see that the probability
that there exists any pair of vertices with an extension to $F$ is
$o(n^{-\gamma})$.  Since $C$ is a constant, the number of
non-isomorphic ways to form $F$ (a superposition of $C$ distinct
copies of $H_1$, overlapping on one particular edge) is still a
constant.  Taking another union bound over all such $F$, we
complete the proof.  \hfill $\Box$

\begin{corollary}
  \label{cor:hme-const-bound-extensions}
  Suppose that
  $n \ll m \ll n^2$, and let $p = 2m/n^2$.  Let $\hme{H}{k}$ have the
  balanced extension property, and let $v$ and $e$ be the numbers of
  vertices and edges in $\hme{H}{k}$.  Suppose that $n^{v-2} p^e$ is a
  negative power of $n$.  Let us consider $G_{jm}$, where $j$ is an
  arbitrary integer constant.  Then, for any constant $\gamma > 0$,
  there exists a constant $C$ such that with probability $1 -
  o(n^{-\gamma})$, every pair of distinct vertices $\{a, b\}$ of
  $G_{jm}$ completes at most $C$ copies of $\hme{H}{(k-1)}$.
\end{corollary}

\noindent {\bf Proof.}\, Fix a pair $\{a, b\}$.  When counting the
number of copies of $\hme{H}{(k-1)}$ completed by that pair, each copy
arises from an extension pair $(H_1, H_2)$ and an extension of
$H_1$ to $H_2$ at the pair.  In fact, this correspondence is
bijective.  The balanced extension property guarantees that all such
pairs are balanced.  Since $H$ is a fixed graph, only a constant
number of non-isomorphic pairs $(H_1, H_2)$ can arise in this way, so
repeated application of Theorem \ref{thm:balanced-extension}(ii)
completes the proof. \hfill $\Box$

\section{Warm-up}
\label{sec:warm-up}

The purpose of this section is to illustrate on a concrete example the
main ideas and techniques that we will use in our proofs.  We
investigate the first nontrivial case, which is the problem of
avoiding $K_4$ in the Achlioptas process with parameter 2.  This turns
out to be the model for the general case.

\begin{theorem*}
  The threshold for avoiding $K_4$ in the Achlioptas process with parameter 2 is $n^{28/19}$.
\end{theorem*}

\noindent {\bf Proof.}\,\, {\em Lower bound:}\, We need to specify a
strategy, and prove that it avoids $K_4$ for many rounds.  At any
intermediate stage in the process, consider a pair of points to be
\emph{2-dangerous}\/ if the addition of an edge between them will
create a copy of $K_4$.  Otherwise, if the addition of the edge will
create a copy of $\hme{K_4}{}$, call the pair \emph{1-dangerous}.
Every other pair is considered to be \emph{0-dangerous}\/ (not
dangerous).  The strategy is then to make an arbitrary choice among
the incoming edges that are minimally dangerous.

Let $m$ be a function of $n$ that satisfies $m \ll n^{28/19}$.  
It suffices to show that for any such $m$, this strategy succeeds \whp.  We
also may assume without loss of generality that $m \gg n^{28/19} / \log n$.  The precise form
of the lower bound on $m$ is not essential; it simplifies the argument by disposing of
uninteresting pathological cases when $m$ is too small.  As it is easier to work
with $G(n, p)$, we will make all of our computations with respect to $p$, which we define
to be $2m/n^2$.  Note that $n^{-10/19}/\log n \ll p \ll n^{-10/19}$.
The following three claims analyze the performance of our strategy.
\begin{description}
\item[(i)] With probability $1 - o(n^{-4})$, $G_m$ has $O(n^4 p^4)$ copies of
  $\hme{K_4}{2}$ and every pair of vertices completes $O(1)$ copies of
  $\hme{K_4}{}$.
\item[(ii)] With probability $1 - o(n^{-2})$, $G_m$ has $O(n^6 p^9)$ 
  copies of $\hme{K_4}{}$.
\item[(iii)] The probability of failure in $m$ rounds is $o(1)$.
\end{description}

For (i), it is easy to verify that $\hme{K_4}{2}$ is a balanced graph,
no matter which two edges are deleted.  Then the number of copies of
$\hme{K_4}{2}$ is roughly what it should be in the random graph $G(n,
p)$---this is made precise by Theorem \ref{thm:balanced-graph}, which
bounds the number of copies of $\hme{K_4}{2}$ in $G_m$ by $O(n^4
p^4)$ \wep\ since $n^4 p^4$ is a positive power of $n$.  It is also
easy to verify that $\hme{K_4}{2}$ has the balanced extension property,
so since $n^2 p^4$ is a negative
power of $n$, Corollary \ref{cor:hme-const-bound-extensions} shows
that there is some constant $C$ such that with probability $1 -
o(n^{-4})$, every pair of vertices in $G_m$ completes at most $C$ copies
of $\hme{K_4}{}$.  This proves (i).

For (ii), fix some $i < m$ and consider the $(i+1)$-st round.  In this
round, the strategy will create one or more copies of $\hme{K_4}{}$
only if both incoming edges span pairs that are 1- or 2-dangerous.
The number of such pairs is at most $O(1)$ times the number of copies of
$\hme{K_4}{2}$.  
Since $G_i \subset G_m$, claim (i) shows that with probability $1 -
o(n^{-4})$, $G_i$ has $O(n^4 p^4)$ copies of $\hme{K_4}{2}$ and
every pair of vertices completes $O(1)$ copies of $\hme{K_4}{}$.
Call this event $A_i$, and condition on it.
Even after conditioning, the incoming edges
at the $(i+1)$-st round are still independently and uniformly 
distributed over the $\Omega(n^2)$ unoccupied pairs of $G_i$, so
the probability that
we are forced to create a new copy of $\hme{K_4}{}$ in this round is
$O\big( \big(\frac{n^4 p^4}{n^2}\big)^2 \big) = O(n^4p^8)$.  
Furthermore, each time this occurs, we only create $O(1)$ new copies 
of $\hme{K_4}{}$ because of our conditioning.  Therefore, the number of new
copies of $\hme{K_4}{}$ in the $(i+1)$-st round is stochastically
dominated by $O(1)$ times the Bernoulli random variable with
parameter $O(n^4p^8)$.  Letting $i$ run through all $m$ rounds, we
see that with probability at least $1 - \sum \pr{\neg A_i} \geq 1 -
o(n^{-2})$, the number of copies of $\hme{K_4}{}$ in $G_m$ is $O(1) \cdot
\bin\big[m, O(n^4p^8)\big]$.  Since $m=n^2p/2$, the expectation of this binomial is a
positive power of $n$, so the Chernoff bound implies that \wep, it is $O(m \cdot n^4p^8) = O(n^6 p^9)$.
This proves (ii).

For (iii), fix some $i$ and consider the probability that we lose in
the $(i+1)$-st round.  The strategy fails precisely when both of the
incoming edges span pairs that are 2-dangerous (completing $K_4$), and
the number of such pairs is at most $O(1)$ times the number of copies of $\hme{K_4}{}$.
Since $G_i \subset G_m$, claim (ii) shows that with probability $1 -
o(n^{-2})$, $G_i$ has $O(n^6 p^9)$ copies of $\hme{K_4}{}$.  Call 
this event $B_i$, and condition on it.
Even after conditioning,
the incoming edges are still independently and uniformly distributed
over the $\Omega(n^2)$ unoccupied pairs of $G_i$, so the probability that both incoming edges are
2-dangerous is $O\big( \big(\frac{n^6 p^9}{n^2}\big)^2 \big) = O(n^8p^{18})$.
Therefore, letting $i$ run through all $m=n^2p/2$ rounds, a union
bound shows that the probability that we are forced to complete a copy
of $K_4$ by the end of the $m$-th round is $\mathbb{P} \leq O(n^2
p \cdot n^8p^{18}) + \sum \pr{\neg B_i} = O(n^{10} p^{19}) + o(1) =
o(1)$.

\vspace{4mm}

\noindent {\em Upper bound:}\, Now suppose that $m \gg n^{28/19}$.  It suffices to
show that we will lose within the first $4m$ rounds \whp.
Again, we may assume without loss of generality
that $m \ll n^{28/19} \log n$, and we will work in terms of $G(n, p)$ with $p =
2m/n^2$.  Note that $n^{-10/19} \ll p \ll n^{-10/19} \log n$.  Let us specify a
sequence of graphs such that each graph is obtained from the previous one by
adding a single edge: let $H_0 = P_4$ (4-vertex path), $H_1 = C_4$
(4-cycle), $H_2 = \hme{K_4}{}$, and $H_3 = K_4$.  It is easy to verify
that the corresponding pairs $(H_0, H_1)$, $(H_1, H_2)$, and $(H_2, H_3)$ are
all balanced extension pairs.  Our result follows from the following four
claims:
\begin{description}
\item[(i)] $G_m$ always contains $\Omega(n^4 p^3)$ copies of $H_0$.  Also, \wep, every pair
  of vertices in $G_{2m}$ extends $O(n^2 p^3)$ copies of $H_0$ into $H_1$.
\item[(ii)] $G_{2m}$ contains $\Omega(n^4 p^4)$ copies of $H_1$ \whp, and with
  probability $1 - o(n^{-2})$, every pair of vertices in $G_{3m}$ extends $O(1)$ copies of $H_1$ into $H_2$.
\item[(iii)] $G_{3m}$ contains $\Omega(n^6 p^9)$
  copies of $H_2$ \whp, and with probability $1 - o(n^{-2})$, every pair of vertices in
  $G_{4m}$ extends $O(1)$ copies of $H_2$ into $H_3$.
\item[(iv)] The probability of survival through $4m$ rounds is $o(1)$.
\end{description}

\vspace{2mm}

\noindent {\bf Proof of (i).}\,  Since the average degree in $G_m$ is
precisely $2m/n = np \gg 1$, from Lemma \ref{lem:count-paths} we conclude that the number of
4-vertex paths is $\Omega(n (np)^3)$. The second part of this claim follows from Theorem
\ref{thm:balanced-extension}(i) since $(H_0, H_1)$ is a balanced
extension pair and $n^2 p^3$ is a positive
power of $n$. \hfill $\Box$

\vspace{2mm}

\noindent {\bf Proof of (ii).}\,
The second part of (ii) follows from Theorem
\ref{thm:balanced-extension}(ii) since $(H_1, H_2)$ is balanced and
$n^2 p^4$ is a negative power of $n$.
To prove the first part of (ii), consider the
$(i+1)$-st round, where $m \leq i < 2m$.  Regardless of the choice of
strategy, if both incoming edges span pairs that extend $\Omega(n^2 p^3)$
copies of $H_0$ into $H_1$, we will be forced to create $\Omega(n^2 p^3)$ new copies of $H_1$.

By (i), the total number of copies of $H_0$ in
$G_i \supset G_m$ is $\Omega(n^4 p^3)$.  For a pair of vertices $\{a,
b\}$, let $n_{a, b}$ be the number of copies of $H_0$ that $\{a, b\}$
extends to $H_1$.  Recall that this definition does not depend
on the presence of an edge between $a$ and $b$.  Since $G_i \subset G_{2m}$, 
claim (i) shows that \wep, in $G_i$ every $n_{a, b} = O(n^2 p^3)$.  
Call this event $A_i$, and condition on it.

Let us estimate the average value of $n_{a, b}$ over all pairs.
Since $H_0$ differs from $H_1$ at exactly one edge, each copy of $H_0$
has a pair at which it contributes $+1$ to the sum $\sum n_{a, b}$.
Therefore, averaging over all $n \choose 2$ pairs of vertices, we
obtain that the average number of copies of $H_0$ that are extended to
$H_1$ at a pair is $\Omega(n^2 p^3)$.  On the other hand, every pair
of vertices in $G_i$ extends $O(n^2 p^3)$ copies of $H_0$ into
$H_1$.  Therefore, at least a constant fraction $\gamma$ 
(where $\gamma = \Omega(1)$ can be chosen to be the same for all $i$) of
all ${n \choose 2}$ pairs have the property of extending $\Omega(n^2
p^3)$ copies of $H_0$ into $H_1$.  Let $P$ be the set of all such
pairs.  Regardless of the choice of strategy, if both incoming edges
span pairs in $P$, we will be forced to create $\Omega(n^2 p^3)$ copies
of $H_1$.  Since $i=o(n^2) = o(|P|)$ and incoming edges are uniformly
distributed over the ${n \choose 2} - i = (1-o(1)){n \choose 2}$
unoccupied pairs, we conclude that the probability that both
incoming edges span pairs in $P$ is $q \geq (1+o(1)) \gamma^2 =
\Omega(1)$.

Let $i$ run from $m$ to $2m$.  Then, up to an error probability of at
most $\sum \pr{\neg A_i} = o(1)$, the number of copies of $H_1$ in $G_{2m}$ is 
at least $\bin(m, q) \cdot \Omega(n^2 p^3)$.  By the Chernoff bound, the
binomial factor exceeds $mq/2 = \Omega(n^2 p)$ \wep; thus, \whp\, $G_{2m}$
has $\Omega(n^2 p \cdot n^2 p^3) = \Omega(n^4p^4)$ copies of $H_1$.  \hfill $\Box$

\vspace{2mm}

\noindent {\bf Proof of (iii).}\, The second part of (iii) follows from Theorem
\ref{thm:balanced-extension}(ii) since $(H_2, H_3)$ is balanced and $n^2 p^5$ 
is a negative power of $n$.
For the first part of (iii), let us consider the $(i+1)$-st round, with $2m \leq i < 3m$.
Regardless of the choice of strategy, if both
incoming edges span pairs that extend copies of $H_1$ into $H_2$, we will create a
copy of $H_2$.  Let $P$ be the set of all such pairs.  We need a lower bound on $|P|$.
Condition on the event $B$ that $G_{2m}$ contains $\Omega(n^4 p^4)$
copies of $H_1$, which occurs \whp\ by (ii).  Also by (ii), with
probability $1 - o(n^{-2})$, every pair of vertices in $G_i$ only
extends $O(1)$ copies of $H_1$ into $H_2$, since $G_i \subset G_{3m}$.
Call this event $C_i$, and condition on it.

Note that every copy of $H_1$ contributes a pair to $P$ which extends $H_1$ into $H_2$,
namely the pair at which it is missing an edge compared to $H_2$. On the other
hand, every such pair was only counted $O(1)$ times, since every pair 
in $G_i$ extends $O(1)$ copies of $H_1$ into $H_2$. This implies that $|P| = \Omega(n^4 p^4)$. 
The incoming edges are uniformly distributed over all unoccupied pairs.  If at least half of the pairs in $P$
were occupied, then we would have $\Omega(n^4 p^4) \gg n^6 p^9$ copies of $H_2$, which would
already give the conclusion of (iii).  Otherwise, the probability that both
incoming edges span pairs in $P$ (hence forcing the creation of a new copy of $H_2$) is $q \geq
(1+o(1)) \big(\frac{|P|/2}{n^2/2}\big)^2 = \Omega\big( \big(\frac{n^4 p^4}{n^2}\big)^2 \big) = \Omega(n^4p^8)$.

Letting $i$ run from $2m$ to $3m$, we see that with error probability at most $\pr{\neg B}
+ \sum \pr{\neg C_i} = o(1)$, either we already obtained the conclusion of (iii), or
the total number of copies of $H_2$ is at least $\bin(m, q)$.  The expectation of this
binomial is $(n^2p/2)q = \Omega(n^6 p^9)$, which is a positive power of $n$. Hence, by the Chernoff bound,
$G_{3m}$ has $\Omega(n^6 p^9)$ copies of
$H_2$ \whp. \hfill $\Box$

\vspace{2mm}

\noindent {\bf Proof of (iv).}\,  Consider the $(i+1)$-st round, where $3m \leq i < 4m$.  Regardless of
the choice of strategy, if both incoming edges span pairs that
complete copies of $H_3 = K_4$, we lose.  We can lower bound the number of such pairs by $\Omega(n^6 p^9)$ by conditioning on the following events.
Let $D$ be the event that $G_{3m}$ contains $\Omega(n^6 p^9)$
copies of $H_2$, which occurs \whp\ by (iii).  Also by (iii), with
probability $1 - o(n^{-2})$, every pair of vertices in $G_i$ extends $O(1)$ copies of
$H_2$ into $H_3$; call this event $E_i$.

Even after conditioning, incoming edges in the $(i+1)$-st round are independently and 
uniformly distributed over the ${n \choose 2} -
i = \Theta(n^2)$ unoccupied pairs of $G_i$.  Therefore, the probability that
both pairs complete $K_4$, conditioned on survival through the
$i$-th round, is $p_i = \Omega\big( \big(\frac{n^6 p^9}{n^2}\big)^2 \big) = \Omega(n^8p^{18})$.  Letting
$i$ run from $3m$ to $4m$, we see that the probability that any strategy
can survive for $4m$ rounds is at most
\begin{eqnarray*}
  \mathbb{P} &\leq& \pr{\neg D} + \sum \pr{\neg E_i} + \prod (1-p_i) 
  \leq o(1)
  + \exp\left\{-\sum p_i\right\} \\
  &\leq& o(1) + \exp\Big\{ -\Omega(n^2 p \cdot n^8p^{18}) \Big\} =
  o(1) + e^{-\omega(1)} \,= \, o(1),
\end{eqnarray*}
which completes the proof. \hfill $\Box$

\section{Avoiding $\boldsymbol{K_t}$, general case}
\label{sec:avoid-kt}

The previous section proved the threshold for avoiding $K_t$ in the
Achlioptas process with parameter $r$, when $t=4$ and $r=2$.  The case
$t=3$ will be covered in Section \ref{sec:avoid-cycles}, which
considers all cycles $C_t$.  In this section, we resolve all other
cases, except for the special case $(t,r) = (4,3)$ which requires more
delicate analysis.
We postpone this final case to Section \ref{sec:avoid,t=4,r=3}.

\begin{theorem*}
  For either $t \geq 5$ and $r \geq 2$, or $t = 4$ and $r \geq 4$, the
  threshold for avoiding $K_t$ in the Achlioptas process with
  parameter $r \geq 2$ is $n^{2-\theta}$, where $\theta$ is defined as
  follows:
  \begin{displaymath}
    s = \lfloor \log_r [(r-1)t + 1] \rfloor, \quad\quad\quad
    \theta = \frac{r^s (t-2) + 2}{r^s\left({t \choose 2} -
    s\right) + \frac{r^s - 1}{r-1}}.
  \end{displaymath}
\end{theorem*}

Before we begin the proof, let us prove an inequality that we will use
in two claims in the lower bound, and the last claim of the upper bound.

\begin{inequality}
  \label{ineq:xy-seq-pos-power}
  Let $a > 2$, $b > 0$, and $r > 1$, and let $s$ be a positive
  integer.  Define the sequences $\{x_s, x_{s-1}, \ldots, x_0\}$ and
  $\{y_s, y_{s-1}, \ldots, y_0\}$ as follows.  Set $x_s = a$ and $y_s
  = b$, and define the rest of the terms recursively by
  \begin{displaymath}
    x_{k-1} = 2 + (x_k - 2)r, \quad \quad \quad y_{k-1} = 1 + y_k r.
  \end{displaymath}
  Then for any $p \gg n^{-x_0/y_0}$, $n^{x_k} p^{y_k}$ is a
  positive power of $n$ for every $k \in \{s, \ldots, 1\}$.
\end{inequality}

\noindent {\bf Proof.}\, Fix any $k \in \{s, \ldots, 1\}$.  One can
easily solve the recursions for $x_k$ and $y_k$ to find:
\begin{displaymath}
  x_k = r^{s-k} (a - 2) + 2, \quad \quad \quad y_k = r^{s-k} b + \frac{r^{s-k} - 1}{r-1}.
\end{displaymath}
Therefore,
\begin{displaymath}
  \frac{x_k}{y_k} = \frac{r^{s-k} (a - 2) + 2}{r^{s-k} b + \frac{r^{s-k} - 1}{r-1}} 
  = \frac{r^s (a - 2) + 2 r^k}{r^s b + \frac{r^s - r^k}{r-1}}.
\end{displaymath}
By the original definition via the recursions, $x_k$ and $y_k$ are
both positive, so the numerator and denominator of the final fraction
above are positive.  Yet as $k$ decreases, the numerator decreases and
the denominator increases.  Therefore, $x_k/y_k > x_0/y_0$.  In
particular, since we assumed that $p \gg n^{-x_0/y_0}$, we conclude
that $n^{x_k} p^{y_k}$ is a positive power of $n$, as desired. \hfill
$\Box$

\vspace{4mm}

Note that if we choose $a = v(K_t) = t$ and $b = e(K_t) - s = {t \choose 2} - s$, then the
above recursions produce $x_0$ and $y_0$ such that the fraction
$x_0/y_0$ is equal to our $\theta$.  Let us now return to the proof of
our thresholds for avoiding $K_t$.

\vspace{4mm}

\noindent {\bf Proof of Theorem.}\,\, {\em Lower bound:}\, The strategy is a natural extension of the one
used to avoid $K_4$.  At any intermediate stage in the process, for any $1 \leq d \leq s$,
consider a pair of points to be \emph{$d$-dangerous}\/ if $d$ is the maximal integer such
that the addition of an edge between them will create a copy of $\hme{K_t}{(s-d)}$.  If
there is no such $d$, consider the pair to be \emph{0-dangerous}.  The
strategy is then to make an arbitrary choice among the incoming edges that are
minimally dangerous.

Let $m \ll n^{2-\theta}$, and let $p = 2m/n^2$.  Again, we assume without loss of generality
that $m \gg n^{2-\theta}/\log n$.  Note that $n^{-\theta}/\log n \ll p \ll
n^{-\theta}$.  We will analyze the performance of our strategy by proving three successive
claims:
\begin{description}
\item[(i)] With probability $1 - o(n^{-2s})$, $G_m$ has $O\big( n^t p^{{t \choose 2} - s} \big)$ copies of
  $\hme{K_t}{s}$, and every pair of vertices completes $O(1)$ copies of
  $\hme{K_t}{(s-1)}$.
\item[(ii)] For each $k \in \{s, s-1, \ldots, 2\}$, and constants
  $x$ and $y$ such that $(n^2 p)\big(\frac{n^x p^y}{n^2}\big)^r$ is a
  positive power of $n$, statement (a) implies statement (b), which
  are defined as follows:
  \begin{description}
  \item[(a)] With probability $1 - o(n^{-2k})$, $G_m$ has $O(n^x
    p^y)$ copies of $\hme{K_t}{k}$, and every pair of vertices
    completes $O(1)$ copies of $\hme{K_t}{(k-1)}$.
  \item[(b)] With probability $1 - o(n^{-2(k-1)})$, $G_m$ has $O\big(
    (n^2 p)\big(\frac{n^x p^y}{n^2}\big)^r \big)$ copies of
    $\hme{K_t}{(k-1)}$, and every pair of vertices completes $O(1)$
    copies of $\hme{K_t}{(k-2)}$.
  \end{description}
\item[(iii)] The probability of failure in $m$ rounds is $o(1)$.
\end{description}
Again, we separate the proofs of the claims for clarity.  At several
points, we require certain inequalities whose rather tedious proofs 
would interfere with the exposition.
\shortversion{%
  The appendix contains the precise formulations of these statements.
}%
{%
  The appendix contains the precise formulations and proofs of these
  statements.
}

\vspace{2mm}

\noindent {\bf Proof of (i).}\, Lemma
\ref{lem:kt-lower-i-balanced-graph} verifies that $\hme{K_t}{s}$ is a
balanced graph, and the $k=s$ case of Inequality
\ref{ineq:xy-seq-pos-power} shows that $n^t p^{{t \choose 2} - s}$
is a positive power of $n$, so Theorem \ref{thm:balanced-graph} implies
that the number of copies of $\hme{K_t}{s}$ in $G_m$ is $O\big(
n^t p^{{t \choose 2} - s} \big)$ \wep.
For the second part of claim (i), Lemma
\ref{lem:kt-lower-i-balanced-extension} verifies that $\hme{K_t}{s}$
has the balanced extension property, and Inequality 
\ref{ineq:kt-upper-neg-power} shows that $n^{t-2} p^{{t \choose 2} - s}$
is a negative power of $n$.  So,
Corollary \ref{cor:hme-const-bound-extensions} shows that there is
some constant $C$ such that with probability $1 - o(n^{-2s})$, every
pair of vertices in $G_m$ completes at most $C$ copies of
$\hme{K_t}{(s-1)}$.  This finishes claim (i). \hfill $\Box$

\vspace{2mm}

\noindent {\bf Proof of (ii).}\, Fix $k$, $x$, and $y$ as specified, and let
us show that (a) implies (b).  First, since every graph of the form
$\hme{K_t}{(k-2)}$ always contains some graph of the form
$\hme{K_t}{(k-1)}$, (a) immediately implies that with probability $1 -
o(n^{-2k})$, every pair of vertices completes $O(1)$ copies of
$\hme{K_t}{(k-2)}$; this implies the second part of (b).

It remains to show the first part of (b).  Fix some $i < m$ and
consider the $(i+1)$-st round.  In this round, the strategy will
create one or more copies of $\hme{K_t}{(k-1)}$ only if all $r$
incoming edges span pairs that are at least $(s-k+1)$-dangerous (i.e.,
create copies of $\hme{K_t}{(k-1)}$).  The number of such pairs is
at most $O(1)$ times the number of copies $\hme{K_t}{k}$.  
Since $G_i \subset G_m$, statement (a) implies that with probability
$1-o(n^{-2k})$, $G_i$ has $O(n^x p^y)$ copies of $\hme{K_t}{k}$ and
every pair of vertices completes $O(1)$ copies of
$\hme{K_t}{(k-1)}$.  Call this event $A_i$, and condition on it.
Even
after conditioning, incoming edges are still independently and
uniformly distributed over the $\Omega(n^2)$ unoccupied pairs of $G_i$, so the
probability that some new copies of $\hme{K_t}{(k-1)}$ are created in
this round is $O\big( \big(\frac{n^x p^y}{n^2}\big)^r \big)$.  Also, by our
conditioning, the number of newly created copies of $\hme{K_t}{(k-1)}$
is still $O(1)$ even when this occurs.  Therefore, the number of
new copies of $\hme{K_t}{(k-1)}$ in the $(i+1)$-st round is
stochastically dominated by $O(1)$ times the Bernoulli random
variable with parameter $O\big( \big(\frac{n^x p^y}{n^2}\big)^r \big)$.  Letting
$i$ run through all $m$ rounds, we see that with probability at least
$1 - \sum \pr{\neg A_i} \geq 1 - o(n^{-2(k-1)})$, the number of copies
of $\hme{K_t}{(k-1)}$ in $G_m$ is $O(1) \cdot \bin\big[m, O\big( \big(\frac{n^x
  p^y}{n^2}\big)^r\big) \big]$.  Since this binomial has expectation $\frac{n^2 p}{2} \cdot O\big( \big(\frac{n^x p^y}{n^2}\big)^r \big)$, which is a positive power of $n$
by the assumption on $x$ and $y$, a Chernoff bound implies that it is
$O\big( (n^2 p)\big(\frac{n^x p^y}{n^2}\big)^r \big)$ \wep.  This finishes
(ii). \hfill $\Box$

\vspace{2mm}

\noindent {\bf Proof of (iii).}\, The idea is to apply claim (i), and then
to repeatedly apply claim (ii) until we obtain a high-probability upper
bound on the number of copies of $\hme{K_t}{}$.  Then, we complete the
proof with essentially the same argument as in claim (iii) of the
proof of the lower bound for avoiding $K_4$.

To keep track of the exponents of $n$ and $p$ in the successive upper
bounds, define the sequences $\{x_s, x_{s-1}, \ldots, x_0\}$ and
$\{y_s, y_{s-1}, \ldots, y_0\}$ as in Inequality \ref{ineq:xy-seq-pos-power}, which then verifies that
$n^{x_k} p^{y_k}$ is a positive power of $n$ for every $k \in \{s-1,
\ldots, 1\}$.  Hence we can apply claims (i) and (ii) until we conclude that
with probability $1 - o(n^{-2})$, $G_m$ has $O(n^{x_1} p^{y_1})$
copies of $\hme{K_t}{}$.  

Now fix some $i$ and consider the probability that we lose in
the $(i+1)$-st round.  The strategy fails precisely when all $r$ of
the incoming edges span pairs that are $s$-dangerous (completing
$K_t$), and the number of such pairs is at most $O(1)$ times the number of copies of
$\hme{K_t}{}$.  
Yet since $G_i \subset G_m$, the previous paragraph shows that with
probability $1-o(n^{-2})$, $G_i$ has $O(n^{x_1} p^{y_1})$ copies of
$\hme{K_t}{}$.  Call this event $B_i$, and condition on it.
Even after
conditioning, incoming edges are still independently and uniformly
distributed over the $\Omega(n^2)$ unoccupied pairs of $G_i$, so the probability
that all incoming edges complete $K_t$ is $O\big( \big(\frac{n^{x_1}
  p^{y_1}}{n^2}\big)^r \big)$.  Therefore, letting $i$ run through all $m =
n^2 p / 2$ rounds, a union bound shows that the probability that 
we are forced to complete a copy of $K_t$ is $\mathbb{P} \leq O\big( (n^2
p)\big(\frac{n^{x_1} p^{y_1}}{n^2}\big)^r \big) + \sum \pr{\neg B_i} =
O(n^{x_0} p^{y_0}) + o(1)$.  This in turn is $o(1)$ because we assumed
that $p \ll n^{-\theta}$ with $\theta = x_0/y_0$.  This completes the
proof. \hfill $\Box$

\vspace{4mm}

\noindent {\em Upper bound:}\, Let $m \gg n^{2-\theta}$, and let
$p = 2m/n^2$.  We will show that \whp, any strategy fails within $\Theta(m)$
rounds, which we again break into periods of length $m$.  We may
assume that $m \ll n^{2-\theta} \log n$ without loss of generality.  Note
that $n^{-\theta} \ll p \ll n^{-\theta} \log n$.  

As in the proof of the upper bound for avoiding $K_4$, we will specify
a sequence of graphs such that each graph is obtained from the previous one by
adding a single edge.  Let $H_1 = K_{\lfloor \frac{t}{2} \rfloor, \lceil
  \frac{t}{2} \rceil}$ (the largest bipartite subgraph of $K_t$), and
arbitrarily choose the rest of the sequence $\{H_2, H_3, \ldots,
H_f\}$, where $H_f = K_t$, by adding one missing edge at a time.
So, $f = 1+{t \choose 2} - \lfloor \frac{t}{2} \rfloor \lceil \frac{t}{2}
\rceil$, which is a constant because we assumed $t$ to be fixed.
Our result follows from the following five claims:
\begin{description}
\item[(i)] $G_m$ contains $\Omega(n^t p^{e(H_1)})$ copies of $H_1$ \whp.

\item[(ii)] Let $k$ be a positive integer for which $n^{t-2}
  p^{e(H_{k-1})}$ is a positive power of $n$.  Then $G_{km}$
  contains $\Omega(n^t p^{e(H_k)})$ copies of $H_k$ \whp.

\item[(iii)] $G_{(f-s)m}$ contains $\Omega(n^t p^{e(H_{f-s})})$ copies
  of $H_{f-s}$ \whp.  Also, $n^{t-2} p^{e(H_{f-s})}$ is a negative
  power of $n$; hence
  with probability $1 - o(n^{-2})$, every pair of vertices in
  $G_{(f-s+1)m}$ extends $O(1)$ copies of $H_{f-s}$ into
  $H_{f-s+1}$.

\item[(iv)] For each $k \in \{s, s-1, \ldots, 2\}$, and constants $x$
  and $y$ such that $n^x p^y \ll n^2$ and $(n^2 p)\big(\frac{n^x
    p^y}{n^2}\big)^r$ is a positive power of $n$, statement (a)
  implies statement (b), which are defined as follows:
  \begin{description}
  \item[(a)] $G_{(f-k)m}$ contains $\Omega(n^x p^y)$ copies of $H_{f-k}$
    \whp, and with probability $1 - o(n^{-2})$, every pair of vertices
    in $G_{(f-k+1)m}$ extends $O(1)$ copies of $H_{f-k}$ into
    $H_{f-k+1}$.

  \item[(b)] $G_{(f-k+1)m}$ contains $\Omega\big( (n^2 p)\big(\frac{n^x
      p^y}{n^2}\big)^r \big)$ copies of $H_{f-k+1}$ \whp, and with probability
    $1 - o(n^{-2})$, every pair of vertices in $G_{(f-k+2)m}$ extends
    $O(1)$ copies of $H_{f-k+1}$ into $H_{f-k+2}$.
  \end{description}

\item[(v)] The probability of survival through $fm = \Theta(m)$ rounds is $o(1)$.
\end{description}

\vspace{2mm}

\noindent {\bf Proof of (i).}\, We will actually prove that $G_m$
contains $\Omega(n^t p^{e(H_1)})$ copies of $H_1$ with certainty, not
just \whp.  However, the rest of the claims only require a \whp\
result in claim (i), so we keep it there for the purpose of
generality.

Since we assumed that $p \gg n^{-\theta}$ and Inequality
\ref{ineq:kt-upper-t/2-codegree} bounds $-\theta \geq -\big\lfloor
\frac{t}{2} \big\rfloor^{-1}$, Lemma \ref{lem:extremal-count-ktt}
implies that the number of copies of the complete bipartite graph $H_1
= K_{\lfloor \frac{t}{2} \rfloor, \lceil \frac{t}{2} \rceil}$ in any
$m$-edge graph is $\Omega( n^t p^{e(H_1)} )$. \hfill
$\Box$

\vspace{2mm}

\noindent {\bf Proof of (ii).}\, We proceed inductively.  The base
case of the induction follows from claim (i).  Now, suppose $k$
satisfies the property that $n^{t-2} p^{e(H_{k-1})}$ is a positive
power of $n$, and $G_{(k-1)m}$ contains $\Omega(n^t p^{e(H_{k-1})})$
copies of $H_{k-1}$ \whp.  We will show that $G_{km}$ contains $\Omega(n^t p^{e(H_k)})$ copies of $H_k$ \whp.

Let us begin by conditioning on the high-probability event $A$ from
our inductive assumption: that $G_{(k-1)m}$ contains $\Omega(n^t
p^{e(H_{k-1})})$ copies of $H_{k-1}$.  Now consider the $(i+1)$-st
round, where $(k-1)m \leq i < km$.  Since $G_i \supset G_{(k-1)m}$,
the total number of copies of $H_{k-1}$ in $G_i$ is $\Omega(n^t
p^{e(H_{k-1})})$ by our conditioning.

Lemma \ref{lem:kt-upper-balanced-extension} verifies that
$(H_{k-1}, H_k)$ is a balanced extension pair, and we assumed that
$n^{t-2} p^{e(H_{k-1})}$ was a positive power of $n$, so Theorem
\ref{thm:balanced-extension}(i) establishes that \wep, every pair of
vertices in $G_{km}$ extends $O(n^{t-2} p^{e(H_{k-1})})$ copies of
$H_{k-1}$ into $H_k$.  Since $G_i \subset G_{km}$, the same bound
holds for $G_i$ \wep; call that event $B_i$, and condition on it.

For a pair of vertices $\{a, b\}$, let $n_{a, b}$ be the number of
copies of $H_{k-1}$ that the pair $\{a, b\}$ extends into $H_k$.
Recall that this definition does not depend on the presence of an edge
between $a$ and $b$.  Let us estimate the average value of $n_{a, b}$
over all pairs.  Since $H_{k-1}$ differs from $H_k$ at exactly one
edge, each copy of $H_{k-1}$ has a pair at which it contributes $+1$
to the sum $\sum n_{a, b}$.  Therefore, averaging over all $n \choose
2$ pairs of vertices, we obtain that the average number of copies of
$H_{k-1}$ that are extended to $H_k$ at a pair is $\Omega(n^{t-2}
p^{e(H_{k-1})})$.

On the other hand, every pair of vertices in $G_i$
extends $O(n^{t-2} p^{e(H_{k-1})})$ copies of $H_{k-1}$ into $H_k$.
Therefore, at least a constant fraction $\gamma = \Omega(1)$ of all $n
\choose 2$ pairs have the property of extending $\Omega(n^{t-2}
p^{e(H_{k-1})})$ copies of $H_{k-1}$ into $H_k$.  Let $P$ be the set of
all such pairs.  Regardless of the choice of strategy, if all $r$
incoming edges span pairs in $P$, we will be forced to create $\Omega(n^{t-2} p^{e(H_{k-1})})$ copies of $H_k$.  Since $i = o(n^2) = o(|P|)$
and incoming edges are uniformly distributed over the ${n \choose 2} -
i = (1+o(1)) {n \choose 2}$ unoccupied pairs, we conclude
that the probability that all incoming edges span pairs in $P$ is $q
\geq (1+o(1)) \gamma^r = \Omega(1)$.

Let $i$ run from $(k-1)m$ to $km$.  Then, up to an error probability
of at most $\pr{\neg A} + \sum \pr{\neg B_i} = o(1)$, the number of
copies of $H_k$ in $G_{km}$ is at least $\bin(m, q) \cdot \Omega(n^{t-2}
p^{e(H_{k-1})})$.  By the Chernoff bound, the binomial factor exceeds
$mq/2 = \Omega(n^2 p)$ \wep; thus, \whp\ $G_{km}$ has $\Omega(n^2
p \cdot n^{t-2} p^{e(H_{k-1})}) = \Omega(n^t p^{e(H_k)})$ copies of $H_k$.
\hfill $\Box$

\vspace{2mm}

\noindent {\bf Proof of (iii).}\, The first part follows directly from claim
(ii), because Inequality \ref{ineq:kt-upper-pos-power} verifies that
$n^{t-2} p^{e(H_{(f-s)-1})}$ is a positive power of $n$.  For the
second part, $(H_{f-s}, H_{f-s+1})$ is a balanced extension pair
by Lemma \ref{lem:kt-upper-balanced-extension}, and $n^{t-2} p^{e(H_{f-s})}$ is a negative
power of $n$ by Inequality \ref{ineq:kt-upper-neg-power}.  Therefore,
Theorem \ref{thm:balanced-extension}(ii) shows that there is some
constant $C$ such that with probability $1 - o(n^{-2})$, every pair of
vertices in $G_{(f-s+1)m}$ extends at most $C$ copies of $H_{f-s}$ into
$H_{f-s+1}$.  This finishes claim (iii). \hfill $\Box$

\vspace{2mm}

\noindent {\bf Proof of (iv).}\, Fix $k$, $x$, and $y$ as specified in the
statement, and assume statement (a).  Let us begin by establishing the
second part of (b).  Lemma \ref{lem:kt-upper-balanced-extension}
verifies that $(H_{f-k+1}, H_{f-k+2})$ is a balanced extension
pair, and Inequality \ref{ineq:kt-upper-neg-power} shows that 
$n^{t-2} p^{e(H_{f-k+1})}$ is a negative power of $n$ for $k \leq s$.  Therefore,
Theorem \ref{thm:balanced-extension}(ii) shows that there is some
constant $C$ such that with probability $1 - o(n^{-2})$, every pair of
vertices in $G_{(f-k+2)m}$ extends at most $C$ copies of $H_{f-k+1}$ into
$H_{f-k+2}$.  This finishes the second part of (b).

It remains to prove the first part of (b).  Consider the $(i+1)$-st
round, with $(f-k)m \leq i < (f-k+1)m$.  Regardless of the choice of
strategy, if all $r$ incoming edges span pairs that extend copies of
$H_{f-k}$ into $H_{f-k+1}$, we will create a copy of $H_{f-k+1}$.  Let
$P$ be the set of all such pairs.  We need a lower bound on $|P|$.

Condition on the high-probability event $C$ of (a) that $G_{(f-k)m}$
contains $\Omega(n^x p^y)$ copies of $H_{f-k}$.  Since $G_i \subset
G_{(f-k+1)m}$, (a) implies that with probability $1 - o(n^{-2})$,
every pair of vertices in $G_i$ extends $O(1)$ copies of $H_{f-k}$
into $H_{f-k+1}$.  Call this event $D_i$, and condition on it.

Note that every copy of $H_{f-k}$ contributes a pair to $P$ which
extends $H_{f-k}$ into $H_{f-k+1}$, namely the pair at which it is
missing an edge compared to $H_{f-k+1}$.  On the other hand, every
such pair was counted at most a constant number of times, since every
pair in $G_i$ extends $O(1)$ copies of $H_{f-k}$ into $H_{f-k+1}$.
This implies that $|P| = \Omega(n^x p^y)$.  The incoming edges are
uniformly distributed over all unoccupied pairs.  If at least half of
the pairs in $P$ were occupied, then we would have $\Omega(n^x
p^y)$ copies of $H_{f-k+1}$.  Yet this would already give us the
conclusion of (b) since:
\begin{displaymath}
  n^x p^y \gg (n^2 p)\Big(\frac{n^x p^y}{n^2}\Big) \gg (n^2 p)\Big(\frac{n^x p^y}{n^2}\Big)^r.
\end{displaymath}
(The first inequality is because $p \ll 1$, and the second inequality
follows from the assumption that $n^x p^y \ll n^2$.)  Otherwise, if
less than half of the pairs in $P$ are occupied, then the probability
that all incoming edges span pairs in $P$ (hence forcing the creation
of a copy of $H_{f-k+1}$) is $q \geq (1+o(1))
\big(\frac{|P|/2}{n^2/2}\big)^r = \Omega\big( \big(\frac{n^x
  p^y}{n^2}\big)^r \big)$.

Letting $i$ run from $(f-k)m$ to $(f-k+1)m$, we see that with error
probability at most $\pr{\neg C} + \sum \pr{\neg D_i} = o(1)$, either
we already obtained the conclusion of (b), or the total number of
copies of $H_{f-k+1}$ is at least $\bin(m, q)$.  The expectation of the
binomial is $\big(\frac{n^2 p}{2}\big) q = \Omega\big( (n^2 p)\big(\frac{n^x
  p^y}{n^2}\big)^r \big)$, which is a positive power of $n$ by assumption.
Hence, by the Chernoff bound, $G_{(f-k+1)m}$ has $\Omega\big( (n^2
p)\big(\frac{n^x p^y}{n^2}\big)^r \big)$ copies of $H_{f-k+1}$ \whp. \hfill
$\Box$

\vspace{2mm}

\noindent {\bf Proof of (v).}\, The result of claim (iii) plugs in
directly to claim (iv), which we may iterate until it gives us a a
lower bound on the number of copies of $H_{f-1} = \hme{K_t}{}$ and an
upper bound on the number of copies of $H_{f-1}$ that any pair extends
into $H_f = K_t$.

To keep track of exponents in the successive lower bounds, define
the sequences $\{x_s, x_{s-1}, \ldots, x_0\}$ and $\{y_s, y_{s-1},
\ldots, y_0\}$ exactly as in Inequality \ref{ineq:xy-seq-pos-power}.
To verify that we can indeed iterate claim (iv), we must show that for
all $k \in \{s, s-1, \ldots, 2\}$, we have that $n^{x_k} p^{y_k} \ll
n^2$, and $n^{x_{k-1}} p^{y_{k-1}}$ is a positive power of $n$.  The
first statement follows from an easy induction: claim (iii) establishes it
for $k=s$, and if $n^{x_k} p^{y_k} \ll n^2$, then $\frac{n^{x_k}
  p^{y_k}}{n^2} \ll 1$, so combined with $p \ll 1$, we see that
$n^{x_{k-1}} p^{y_{k-1}} = (n^2 p)\big(\frac{n^{x_k}
  p^{y_k}}{n^2}\big)^r \ll n^2$.  The second statement is verified by
Inequality \ref{ineq:xy-seq-pos-power}.  Therefore, we arrive at
the result that $G_{(f-1)m}$ contains $\Omega(n^{x_1} p^{y_1})$ copies
of $H_{f-1} = \hme{K_t}{}$ \whp.  Call this event $E$, and condition on it.
We also find that with probability $1 - o(n^{-2})$,
every pair of vertices in $G_{fm}$ extends $O(1)$ copies of
$H_{f-1}$ into $H_f$ (i.e., completes $O(1)$ copies of $K_t$).  
The same probability bound also holds in $G_i$ for any $i \leq fm$, 
because $G_i \subset G_{fm}$; let $F_i$ be the corresponding event.

Now consider the $(i+1)$-st round, where $(f-1)m \leq i < fm$.
Regardless of the choice of strategy, if all $r$ incoming edges span
pairs that complete copies of $K_t$, we will lose.  
We can bound the number of such pairs by $\Omega(n^{x_1} p^{y_1})$
by conditioning on the above events $E$ and $F_i$.
Even after conditioning, the incoming edges in this round
still independent and uniformly distributed over the ${n \choose 2} -
i = \Theta(n^2)$ unoccupied pairs of $G_i$.  Therefore, the probability that
all $r$ pairs complete $K_t$, conditioned on survival through the
$i$-th round, is $p_i = \Omega\big( \big(\frac{n^{x_1} p^{y_1}}{n^2}\big)^r \big)$.  Letting
$i$ run from $(f-1)m$ to $fm$, we see that the probability that any strategy 
can survive for $fm$ rounds is at most
\begin{eqnarray*}
  \mathbb{P} &\leq& \pr{\neg E} + \sum \pr{\neg F_i} + \prod (1-p_i) 
  \leq o(1)
  + \exp\left\{-\sum p_i\right\} \\
  &\leq& o(1) + \exp\Big\{ -\Omega\Big((n^2 p)\Big(\frac{n^{x_1} p^{y_1}}{n^2}\Big)^r\Big) \Big\} = o(1) + \exp\{-\Omega(n^{x_0} p^{y_0})\}.
\end{eqnarray*}
This in turn is $o(1)$ because we assumed that $p \gg n^{-\theta}$
with $\theta = x_0/y_0$.  This completes the proof. \hfill $\Box$

\section{Abstraction into general argument}
\label{sec:abstract-argument}

Note that we structured our exposition of the previous section in the
following manner.  The arguments did not directly use properties of the
specific graph that we were avoiding ($K_t$).  Rather, they were linked
to lemmas and inequalities that proved certain properties (e.g.,
balanced-ness, etc.) about $K_t$.  Let us now isolate these necessary
``ingredients'' that one can plug in to our general machinery to prove
thresholds.

For the rest of this section, let $H$ be the fixed graph which we wish to avoid.  Our
arguments allow one to prove the threshold for avoiding $H$ in the
Achlioptas process with parameter $r$ simply by specifying several
parameters, and then proving some lemmas and inequalities that do not need to
refer to the Achlioptas process at all.  We first describe the
parameters.
\begin{itemize}
\item $s$: this was the number of levels of danger considered by the
  avoidance strategy in the proof of the lower bound.  At any
  intermediate stage in the process, for any $1 \leq d \leq s$, we
  considered a pair of points to be \emph{$d$-dangerous}\/ if $d$ was
  the maximal integer such that the addition of an edge between them
  created a copy of $\hme{H}{(s-d)}$.  If there was no such $d$,
  we considered the pair to be \emph{0-dangerous}.
  Recall that the strategy was then to make an arbitrary choice among
  the incoming edges that were minimally dangerous.

\item A sequence of graphs $\{H_1, \ldots, H_f\}$ sharing the same
  vertex set, with each successive graph containing exactly one more
  edge: this was used in the upper bound argument to iteratively prove
  lower bounds on the number of copies of $H_i$, proceeding from $i=1$
  to $i=f$.
\end{itemize}
The correct choice of $s$ then determined $\theta$, the negative
exponent in the threshold (in terms of $p$) for avoidance:
\begin{displaymath}
  \theta = \frac{r^s (v(H)-2) + 2}{r^s\left(e(H) - s\right) + \frac{r^s - 1}{r-1}}.
\end{displaymath}

Assuming that the parameters were suitably chosen, one then only needed
to establish the following lemmas and inequalities in order to prove
that the threshold for avoiding $H$ in the Achlioptas process with
parameter $r$ is $n^{2-\theta}$.

\vspace{2mm}

\noindent\emph{For proof of lower bound.  Here, $n^{-\theta} / \log n \ll p \ll n^{-\theta}$.}
\begin{enumerate}
\item \emph{$\hme{H}{s}$ is a balanced graph.}  This allowed us to
  prove in claim (i) that \wep, $G_m$ has $O(n^{v(H)} p^{e(H)-s})$
  copies of $\hme{H}{s}$.  For $H = K_t$, this was provided by Lemma
  \ref{lem:kt-lower-i-balanced-graph}.

\item \emph{$\hme{H}{s}$ has the balanced extension property, and
    $n^{v(H) - 2} p^{e(H) - s}$ is a negative power of $n$.}  This
  allowed us to prove in claim (i) that with probability $1 -
  o(n^{-2s})$, every pair of vertices in $G_m$ completes $O(1)$
  copies of $\hme{H}{(s-1)}$.  For $H = K_t$, these were provided by Lemma
  \ref{lem:kt-lower-i-balanced-extension} and Inequality
  \ref{ineq:kt-upper-neg-power}.

\end{enumerate}

\noindent\emph{For proof of upper bound.  Here, $n^{-\theta} \ll p \ll n^{-\theta} \log n$.}
\begin{enumerate}
\item \emph{$G_m$ contains $\Omega(n^{v(H_1)} p^{e(H_1)})$ copies of
    $H_1$ \whp.}  This was claim (i), and for $H = K_t$, it was
  provided by the extremal estimate on the number of $K_{s,t}$ (Lemma
  \ref{lem:extremal-count-ktt}), along with Inequality
  \ref{ineq:kt-upper-t/2-codegree}, which assured that $p$ was large
  enough to apply the extremal result.

\item \emph{Each consecutive pair $(H_k, H_{k+1})$ is a balanced
    extension pair.}  This was used throughout the proof of the
  upper bound, and for $H = K_t$, it was provided by Lemma
  \ref{lem:kt-upper-balanced-extension}.

\item \emph{$n^{v(H)-2} p^{e(H)-s-1}$ is a positive power of $n$.}
  This was used in claim (iii) to show that we could iterate the
  argument of claim (ii) enough times to conclude that $G_{(f-s)m}$
  contained $\Omega(n^{v(H_{f-s})} p^{e(H_{f-s})})$ copies of $H_{f-s}$
  \whp.  For $H = K_t$, this was provided by Inequality
  \ref{ineq:kt-upper-pos-power}.

\item \emph{$n^{v(H)-2} p^{e(H)-s}$ is a negative power of $n$.}  This
  was used in claim (iii) to transition to the next inductive process,
  which relied on the copies of $H_{f-s}$ not being too concentrated
  on any pair of vertices.  \textbf{Note:} this statement was already
  required above for the lower bound, so we do not need to
  check it again.
\end{enumerate}

\section{Avoiding cycles}
\label{sec:avoid-cycles}

Now we show by example how to use our machinery to prove avoidance
thresholds.  We start with an easy application which completely solves
the problem for cycles $C_t$.  In light of the previous section, we
only need to provide the required parameters, lemmas, and
inequalities.  We will specify these in the same order that they were
presented in the previous section.  This will prove the following
theorem.

\begin{theorem*}
  For $t \geq 3$, the threshold for avoiding $C_t$ in the Achlioptas
  process with parameter $r \geq 2$ is
  $n^{2-\frac{r(t-2)+2}{r(t-1)+1}}$.
\end{theorem*}

\noindent \textbf{Proof.}\, We use the parameter $s=1$, and the
sequence of graphs $H_1 = \hme{C_t}{}$, $H_2 = C_t$.  This gives the
threshold $n^{2-\theta}$, where $\theta = \frac{r^s (v(C_t)-2) +
  2}{r^s\left(e(C_t) - s\right) + \frac{r^s - 1}{r-1}} =
\frac{r(t-2)+2}{r(t-1)+1}$, which matches the claimed result.  Now we
need to provide the required lemmas and inequalities.  For the
reader's convenience, we have reproduced the italicized statements
from Section \ref{sec:abstract-argument}.

\vspace{2mm}

\noindent\emph{For proof of lower bound.  Here, $n^{-\theta} / \log n \ll p \ll n^{-\theta}$.}
\begin{enumerate}
\item \emph{$\hme{C_t}{}$ is a balanced graph.}  This is obvious.

\item \emph{$\hme{C_t}{}$ has the balanced extension property, and
    $n^{v(C_t) - 2} p^{e(C_t) - 1} = n^{t-2} p^{t-1}$ is a negative
    power of $n$.}  The first part is obvious.  For the second, 
  since $p \ll n^{-\frac{r(t-2)+2}{r(t-1)+1}}$, we must establish that
  $(t-2) - (t-1)\frac{r(t-2)+2}{r(t-1)+1} < 0$.  Routine algebra shows
  that the left hand side equals $-\frac{t}{r(t-1)+1}$, which is
  certainly negative when $t \geq 3$, $r \geq 2$.
\end{enumerate}

\noindent\emph{For proof of upper bound.  Here, $n^{-\theta} \ll p \ll n^{-\theta} \log n$.}
\begin{enumerate}
\item \emph{$G_m$ contains $\Omega(n^{v(H_1)} p^{e(H_1)})$ copies of
    $H_1$ \whp.} The average degree of $G_m$ is precisely $np$ by the
  definition of $p = 2m/n^2$.  We show in item \#3 below that $np$ is
  a positive power of $n$, so it tends to infinity with $n$.  Thus, we
  may apply Lemma \ref{lem:count-paths}, an extremal result counting
  the number of paths, to conclude that $G_m$ contains at least
  $(1+o(1)) n (np)^{t-1}$ copies of the $t$-vertex path $H_1$, as
  desired.

\item \emph{$(H_1, H_2)$ is a balanced extension pair.}  This is
  easy to see.

\item \emph{$n^{v(C_t)-2} p^{e(C_t)-1-1} = (np)^{t-2}$ is a positive
    power of $n$.}  It suffices to show that $np$ is a positive power
  of $n$.  Since $p \gg n^{-\frac{r(t-2)+2}{r(t-1)+1}}$, this amounts
  to proving that $1 - \frac{r(t-2)+2}{r(t-1)+1} > 0$.  Routine
  algebra shows that the left hand side equals $\frac{r-1}{r(t-1)+1}$,
  which is certainly positive when $t \geq 3$, $r \geq 2$.
\end{enumerate}

As we have provided all of the necessary ingredients to apply our
machinery, we are done.  \hfill $\Box$

\section{Avoiding $\boldsymbol{K_{t,t}}$}
\label{sec:avoid-ktt}

Now we show a more complex application of our machinery, which
completely solves the problem for $K_{t,t}$.  This will prove the
following theorem.

\begin{theorem*}
  Suppose that $t \geq 3$ and $r \geq 2$ are fixed integers.
  The threshold for avoiding
  $K_{t,t}$ in the Achlioptas process with parameter $r$ is
  $n^{2-\theta}$, where $\theta$ is defined as follows:
  \begin{displaymath}
    s = \lfloor \log_r [(r-1)t + 1] \rfloor, \quad\quad\quad \theta = \frac{r^s (2t-2) + 2}{r^s (t^2 - s) + \frac{r^s - 1}{r-1}}.
  \end{displaymath}
\end{theorem*}

\subsection{Parameters}

The value of $s$ is already specified in the statement of the theorem,
so we proceed to give the sequence of graphs $\{H_1, \ldots,
H_f\}$.  The sequences are quite different depending on the parity of
$t$, so we describe them separately.

\begin{description}
\item[Case 1: $\boldsymbol{t}$ is even.] Let $H_1$ be the 4-partite graph with
  parts $V_1, V_2, V_3, V_4$, each of size $t/2$, and edges such that
  $(V_1, V_2)$, $(V_1, V_4)$, and $(V_3, V_2)$ are complete
  bipartite graphs.
  Let $\{H_2, \ldots, H_{1+(t/2)}\}$ be obtained by successively
  adding single edges until $H_{1+(t/2)}$ has a perfect matching
  between $V_3$ and $V_4$.  Then, arbitrarily choose the rest of the
  sequence $\{H_{2+(t/2)}, \ldots, H_f\}$ by adding one edge at a
  time, until the final term is the complete bipartite graph $K_{t,t}$
  with bipartition $(V_1 \cup V_3, V_2 \cup V_4)$.  Note that $f =
  1+t^2/4$.

\item[Case 2: $\boldsymbol{t}$ is odd.]  Let $H_1$ be a 6-partite graph with parts
  $\{V_i\}_1^6$ such that $V_3$ and $V_4$ are singletons, and the
  other four parts each have size $\lfloor t/2 \rfloor$.  The
  edges are as follows: 
  the two pairs $(V_1, V_2)$ and $(V_5, V_6)$ are
  each complete bipartite graphs, the vertex in $V_3$ is adjacent
  to all of $V_2 \cup V_4 \cup V_6$, and the vertex in $V_4$ is
  adjacent to all of $V_1 \cup V_3 \cup V_5$.  There are no more
  edges.

  Let $\{H_2, \ldots, H_{1 + \lfloor t/2 \rfloor}\}$ be obtained by
  successively adding single edges until $H_{1 + \lfloor t/2 \rfloor}$
  has a perfect matching between $V_1$ and $V_6$.  To create the next
  $\lfloor t/2 \rfloor$ graphs in the sequence, we put down a matching
  between $V_5$ and $V_2$, one edge at a time.  Finally, arbitrarily
  choose the rest of the sequence $\{H_{2 + 2 \lfloor t/2 \rfloor},
  \ldots, H_f\}$ by adding one edge at a time, until the final term is
  the complete bipartite graph $K_{t,t}$ with bipartition $(V_1 \cup
  V_3 \cup V_5, V_2 \cup V_4 \cup V_6)$.  Note that $f = 1 + 2 \lfloor
  t/2 \rfloor^2$.
\end{description}

\subsection{Lemmas and inequalities}

Next, we provide the required lemmas and inequalities.  For the
reader's convenience, we have reproduced the italicized statements
from Section \ref{sec:abstract-argument}.

\vspace{2mm}

\noindent\emph{For proof of lower bound.  Here, $n^{-\theta} / \log n \ll p \ll n^{-\theta}$.}
\begin{enumerate}
\item \emph{$\hme{K_{t,t}}{s}$ is a balanced graph.}  This is now
  provided by Lemma \ref{lem:ktt-lower-i-balanced-graph}.  Actually,
  the graph is not balanced when $t=3$ and $r=2$, but in that
  particular case, Lemma \ref{lem:ktt-lower-i-balanced-graph}
  additionally proves that the number of copies of $\hme{K_{t,t}}{s}$
  in $G_m$ is still $O(n^{v(H)} p^{e(H)-s})$ \wep, which is all we really need.

\item \emph{$\hme{K_{t,t}}{s}$ has the balanced extension property,
    and $n^{v(K_{t,t}) - 2} p^{e(K_{t,t}) - s}$ is a negative power of
    $n$.}  These are now provided by Lemma
  \ref{lem:ktt-lower-i-balanced-extension} and Inequality 
  \ref{ineq:ktt-upper-neg-power}.
\end{enumerate}

\noindent\emph{For proof of upper bound.  Here, $n^{-\theta} \ll p \ll n^{-\theta} \log n$.}
\begin{enumerate}
\item \emph{$G_m$ contains $\Omega(n^{v(H_1)} p^{e(H_1)})$ copies of
    $H_1$ \whp.}  This time, we use Inequality
  \ref{ineq:ktt-upper-t/2-codegree} to show that $-\theta > -2/t$.
  Since we assume that $p \gg n^{-\theta}$ for the upper bound
  argument, this provides the condition required to apply either Lemma
  \ref{lem:extremal-count-ktt-backbone} if $t$ is even, or Lemma
  \ref{lem:count-ktt-odd-backbone} if $t$ is odd.  Both lemmas (presented below) lead to the required final
  statement.

\item \emph{Each consecutive pair $(H_k, H_{k+1})$ is a balanced
    extension pair.}  This is now provided by Lemma
  \ref{lem:ktt-upper-balanced-extension-even} if $t$ is even, and
  by Lemma 
  \ref{lem:ktt-upper-balanced-extension-odd} if $t$ is odd.

\item \emph{$n^{v(K_{t,t})-2} p^{e(K_{t,t})-s-1}$ is a positive power of $n$.}
  This is provided by Inequality \ref{ineq:ktt-upper-pos-power}.
\end{enumerate}

\subsection{Proofs of supporting lemmas}

We conclude this section by proving the two lemmas that provide the
first component of the proof of the upper bound.  We start with the
lemma that is used when $t$ is even.


\begin{lemma}
  \label{lem:extremal-count-ktt-backbone}
  For any fixed positive integers $k$ and $l$, consider the following
  4-partite graph, which we call $H$.  Let the parts be $V_1, V_2,
  V_3, V_4$, with $|V_1| = |V_2| = k$ and $|V_3| = |V_4| = l$, and
  place edges such that $(V_1, V_2)$, $(V_1, V_4)$, and $(V_2, V_3)$
  are complete bipartite graphs.  There are no more edges.  Then, there
  exists a constant $c_k$ such that for any $p \gg n^{-1/k}$, every
  graph with $n$ vertices and ${n \choose 2} p$ edges contains 
  at least $(c_k+o(1)) n^{2k+2l} p^{k^2 + 2kl}$ copies of $H$.
\end{lemma}

\noindent {\bf Proof.}\, Let us fix an ambient graph $G$ with $n$
vertices and ${n \choose 2} p$ edges.  By Lemma
\ref{lem:extremal-count-ktt}, the number of copies of $K_{k,k}$ in $G$
is at least $(1+o(1)) n^{2k} p^{k^2}$.  Recall that the $k$-codegree of a
set $U$ of $k$ distinct vertices is the number of vertices that are
adjacent to all of $U$.  Let us say that a copy of $K_{k,k}$ is
\emph{deficient}\/ if either of the sides of its bipartition has
$k$-codegree less than $\frac{1}{2} np^k$ in $G$.  We claim that at
most $\frac{1}{2} + o(1)$ of the copies of $K_{k,k}$ are deficient.

To see this, note that if an ordered $k$-tuple of distinct vertices
has $k$-codegree less than $\frac{1}{2} np^k$, then it can extend to at
most $\big(\frac{1}{2} np^k\big)^k$ copies of $K_{k,k}$.  The number
of such $k$-tuples is at most $n^k$; therefore, the number of
deficient copies of $K_{k,k}$ is at most $n^k \big(\frac{1}{2}
np^k\big)^k \leq \frac{1}{2} n^{2k} p^{k^2}$, as claimed.

Yet each non-deficient copy of $K_{k,k}$ extends to at least 
\begin{displaymath}
  {\frac{1}{2}np^k - 2k \choose l}l! \cdot {\frac{1}{2}np^k - 2k - l \choose l}l!
\end{displaymath}
copies of $H$.  This is because we may consider the copy of $K_{k,k}$
to be $V_1 \cup V_2$, we choose $V_3$ from the common neighborhood of
$V_2$ excluding the $2k$ vertices in $V_1 \cup V_2$,
and finally we choose $V_4$ from the common neighborhood of $V_1$
excluding the $2k+l$ vertices in $V_1 \cup V_2 \cup V_3$.
Since we assumed that $p \gg n^{-1/k}$, the binomial coefficients are
asymptotically monomials of degree $l$, so we conclude that each
non-deficient copy of $K_{k,k}$ extends to $\Omega( (np^k)^l \cdot
(np^k)^l ) = \Omega(n^{2l} p^{2kl})$ copies of $H$.  Since there are always
at least $\big(\frac{1}{2} + o(1)\big) n^{2k} p^{k^2}$ non-deficient
copies of $K_{k,k}$, we conclude that the number of copies of $H$ is
always $\Omega(n^{2k+2l} p^{k^2 + 2kl})$, as claimed.  \hfill $\Box$

\vspace{4mm}

Using Lemma \ref{lem:extremal-count-ktt-backbone} as a building block,
we now prove the lemma that provides the first component of the upper
bound when $t$ is odd.  Actually, we prove a result for $G_{2m}$
instead of $G_m$, but this does not matter for the purpose of the
general argument.


\begin{lemma}
  \label{lem:count-ktt-odd-backbone}
  Let $k$ be a positive integer.  Let $H$ be a 6-partite graph with
  parts $\{V_i\}_1^6$ such that $V_3$ and $V_4$ are singletons, and
  the other four parts each have size $k$.  Let
  the edges of $H$ be as follows: the two pairs $(V_1, V_2)$ and $(V_5,
  V_6)$ are each complete bipartite graphs, the vertex in $V_3$ is
  adjacent to all of $V_2 \cup V_4 \cup V_6$, and the vertex in
  $V_4$ is adjacent to all of $V_1 \cup V_3 \cup V_5$.  There are no
  more edges.

  Consider $G_{2m}$, the graph after the $2m$-th round of the
  Achlioptas process with parameter $r \geq 2$.  Let $p = 2m/n^2$, and
  suppose that $p \gg n^{-\theta}$ with $-\theta >
  -1/(k+\frac{1}{2})$.  Then $G_{2m}$ contains $\Omega(n^{v(H)}
  p^{e(H)})$ copies of $H$ \whp.
\end{lemma}

\noindent {\bf Proof.}\, Let $H_1$ be the subgraph of $H$ induced by
$V_1 \cup V_2 \cup V_3 \cup V_4$, and let $H_0$ be the subgraph of
$H_1$ with the edge between $V_3$ and $V_4$ deleted.  Observe that
we can find a copy of $H$ in a graph by first looking for a pair of
vertices for the site of the edge between $V_3$ and $V_4$, and then
looking for two disjoint copies of $H_0$ that are extended into
$H_1$ by that pair.

Consider the $(i+1)$-st turn, for some $m \leq i < 2m$.  By Lemma
\ref{lem:extremal-count-ktt-backbone}, $G_m$ (and hence $G_i \supset
G_m$) always contains $\Omega(n^{2k+2} p^{k^2 + 2k})$ copies of $H_0$.
Lemma \ref{lem:ktt-upper-i-balanced-extension} verifies that
$(H_0, H_1)$ is a balanced extension pair, and $n^{2k} p^{k^2 + 2k}$
is a positive power of $n$ because we assumed that $p \gg n^{-\theta}$
with $-\theta > -1/(k+\frac{1}{2})$ and $k \geq 1$.  Thus, Theorem
\ref{thm:balanced-extension}(i) establishes that \wep, every pair of
vertices in $G_i \subset G_{2m}$ extends $O(n^{2k} p^{k^2 + 2k})$ copies of
$H_0$ into $H_1$.  Call this event $A_i$, and condition on it.

For a pair of vertices $\{a, b\}$, let $n_{a, b}$ be the number of
copies of $H_0$ that the pair $\{a, b\}$ extends into $H_1$.  Recall
that this definition does not depend on the presence of an edge
between $a$ and $b$.  Let us estimate the average value of $n_{a, b}$
over all pairs.  Since $H_0$ differs from $H_1$ at exactly one edge,
each copy of $H_0$ has a pair at which it contributes $+1$ to the sum
$\sum n_{a, b}$.  Therefore, averaging over all $n \choose 2$ pairs of
vertices, we obtain that the average number of copies of $H_0$ that
are extended to $H_1$ at any pair is $\Omega(n^{2k+2} p^{k^2 + 2k})$.

On the other hand, by our conditioning, every pair of vertices in $G_i$ extends $O(n^{2k} p^{k^2 + 2k})$ copies of $H_0$ into $H_1$.  Therefore, at least
a constant fraction $\gamma = \Omega(1)$ of all $n \choose 2$ pairs have
the property of extending $\Omega(n^{2k} p^{k^2 + 2k})$ copies of $H_0$
into $H_1$.  Let $P$ be the set of all such pairs.  Regardless of the
choice of strategy, if all $r$ incoming edges span pairs in $P$, we
will be forced to choose a pair in $P$.  This will create $\Omega\big(
\big(n^{2k} p^{k^2 + 2k}\big)^2 \big) = \Omega(n^{4k} p^{2k^2 + 4k})$ pairs of
copies of $H_0$ that are extended to $H_1$ by the chosen pair.  Such a
pair of copies of $H_0$ would become a new copy of $H$ after the edge
is added, if the pair of copies were disjoint.  If the pair of copies
of $H_0$ is not disjoint, then let us say that they create a
\emph{degenerate}\/ copy of $H$.  For now, let us count degenerate
copies of $H$ along with the true copies of $H$.  Later, we will show
that the degenerate copies are vastly outnumbered by true copies of
$H$.

Since $i = o(n^2) = o(|P|)$ and incoming edges are uniformly
distributed over the ${n \choose 2} - i = (1+o(1)) {n \choose 2}$
unoccupied pairs, we conclude that the probability that all
incoming edges span pairs in $P$ is $q \geq (1+o(1)) \gamma^r = \Omega(1)$.  
Let $i$ run from $m$ to $2m$.  Then \wep, the number of (possibly
degenerate) copies of $H$ in $G_{2m}$ is at least $\bin(m, q)
\cdot \Omega(n^{4k} p^{2k^2 + 4k})$.  By the Chernoff bound, the binomial
factor exceeds $mq/2 = \Omega(n^2 p)$ \wep, so we conclude that $G_{2m}$
has $\Omega(n^2 p \cdot n^{4k} p^{2k^2 + 4k}) = \Omega(n^{v(H)} p^{e(H)})$
(possibly degenerate) copies of $H$ \whp.

To finish the proof of this lemma, we must show that the number
of degenerate copies of $H$ in $G_{2m}$ is $o(n^{v(H)} p^{e(H)})$
\whp.  For convenience, we will work with $G(n, p)$ instead of
$G_{2m}$ because Lemma \ref{lem:coupling} shows that we may couple
$G_{2m}$ with $G(n, 4rp)$, and the constant $4r$ disappears under the
``$o(\cdot)$'' notation.  Note that the underlying graph of a degenerate
copy of $H$ is a superposition of two copies of $K_{k+1,k+1}$, overlapping
on at least 3 vertices.  So, let us consider any such superposition,
and call the underlying graph $F$.  Let $v' = v(F)$ and $e' = e(F)$.
The copies overlap on at least 3 vertices, so $v' < v(H)$.  It is easy to check that since
$K_{k+1,k+1}$ is a balanced graph, $\frac{e'}{v'} \geq \frac{e(H)}{v(H)}$.
So, the expected number of copies of $F$ in $G(n, p)$ is:
\begin{displaymath}
  \mathbb{E} \leq n^{v'} p^{e'} = (n p^{e'/v'})^{v'} \leq (n p^{e(H)/v(H)})^{v'} = (n^{v(H)} p^{e(H)})^{v'/v(H)}.
\end{displaymath}
Now, we assumed that $p \gg n^{-1/(k + \frac{1}{2})}$, so $n^{v(H)}
p^{e(H)} \gg 1$ because $v(H) = 4k+2$ and $e(H) = 2k^2 + 4k + 1$.
Furthermore, $v' < v(H)$, so Markov's inequality implies that
\whp, $G(n, p)$ has $o(n^{v(H)} p^{e(H)})$ copies of $F$.  Since each
copy of $F$ can account for at most a constant number (depending only
on $k$) of degenerate copies of $H$, and there is only a constant
number of non-isomorphic superpositions $F$,
we conclude that \whp, $G(n, p)$ has $o(n^{v(H)} p^{e(H)})$
degenerate copies of $H$.  This completes the proof of the lemma.
\hfill $\Box$

\section{Avoiding $\boldsymbol{K_4}$ in the Achlioptas process with parameter 3}
\label{sec:avoid,t=4,r=3}

To apply the machinery of Section \ref{sec:abstract-argument},
one needs to prove that certain quantities are positive or negative
powers of $n$.  In our study of avoiding cycles, cliques, and complete
bipartite graphs, the only case in which we encounter a key exponent
that is not separated from zero is when we are avoiding $K_4$ in
the Achlioptas process with parameter 3.

However, the separation of the exponent from zero was merely a
convenience which allowed us to bound maxima of families of random
variables (e.g., the maximum codegree in a graph) \whp.  When we do
not have this condition, we may instead bound the entire distribution
of the family.

\begin{lemma}
  \label{lem:codegree-sequence}
  Let $n^{-1/2} \ll p \ll n^{-1/2} \log n$.  Then $G(n, p)$ satisfies
  the following property \whp: all codegrees are at most $np^2 \log
  n$, and for every integer $4 \leq k \leq \log n$, the number of
  pairs with codegree at least $knp^2$ is at most $n^2/k^3$.
\end{lemma}

This result, which we prove at the end of this section, allows us to
prove our final threshold.

\begin{theorem*}
  The threshold for avoiding $K_4$ in the Achlioptas process with
  parameter $3$ is $n^{3/2}$.
\end{theorem*}

\noindent {\bf Proof.}\,\, {\em Lower bound:}\, A shortsighted
strategy works in this instance: arbitrarily select any one of the
incoming edges that does not create a copy of $K_4$.  Let $m \ll
n^{3/2}$, and let $p = 2m/n^2$.  Again, we assume without loss of
generality that $m \gg n^{3/2}/\log n$.  Note that $n^{-1/2}/\log n
\ll p \ll n^{-1/2}$.  We will analyze the performance of our strategy
by proving two successive claims:
\begin{description}
\item[(i)] $G_m$ has $O(n^4 p^5)$ copies of $\hme{K_4}{}$ \wep.
\item[(ii)] The probability of failure in $m$ rounds is $o(1)$.
\end{description}

The interested reader may check that if we followed the recipe for
avoiding $K_t$ in Section \ref{sec:avoid-kt}, we would start by
counting copies of $\hme{K_4}{2}$ instead of $\hme{K_4}{}$.  This is
essentially the only change in the lower bound argument, but we 
provide the details below for completeness.

For (i), $\hme{K_4}{}$ is balanced 
and $n^4 p^5$ is a positive power of $n$, so Theorem
\ref{thm:balanced-graph} implies that the number of copies of
$\hme{K_4}{}$ in $G_m$ is $O(n^4 p^5)$ \wep.

For (ii), consider the probability that the strategy fails at the
$(i+1)$-st round for some $i < m$, i.e., that all 3 incoming edges
span pairs that complete copies of $K_4$.  The number of such pairs is
upper bounded by the number of copies of $\hme{K_4}{}$.  
Since $G_i \subset G_m$, claim (i) implies that $G_i$ has $O(n^4
p^5)$ copies of $\hme{K_4}{}$ \wep.  Call this event $A_i$, and
condition on it.
Then, the chance that all 3 incoming edges
complete $K_4$ is $O\big( \big(\frac{n^4 p^5}{n^2}\big)^3 \big) = O(n^6 p^{15})$.  Letting
$i$ run through all $m = n^2 p/2$ rounds, a union bound shows that the
probability that we are forced to complete a copy of $K_4$ by the $m$-th
round is $\mathbb{P} \leq O(n^2 p \cdot n^6 p^{15}) + \sum \pr{\neg
  A_i} = O(n^8 p^{16}) + o(1) = o(1)$, as desired.

\vspace{4mm}

\noindent {\em Upper bound:}\, Let $m \gg n^{3/2}$, and let $p =
2m/n^2$.  We will show that \whp, any strategy fails within $3m$
rounds, which we break into periods of length $m$.  Again, we may
assume that $m \ll n^{3/2} \log n$ without loss of generality.  Note
that $n^{-1/2} \ll p \ll n^{-1/2} \log n$.  Our result follows from
the following three claims:
\begin{description}
\item[(i)] $G_m$ contains $\Omega(n^2)$ pairs of vertices with codegree
  at least 2 \whp.
\item[(ii)] $G_{2m}$ contains $\Omega(n^2 p)$ copies of $\hme{K_4}{}$
  \whp, and with probability $1 - o(n^{-2})$, every pair of vertices
  in $G_{3m}$ extends $O(1)$ copies of $\hme{K_4}{}$ into $K_4$
\item[(iii)] The probability of survival through $3m$ rounds is $o(1)$.
\end{description}

\vspace{2mm}

\noindent {\bf Proof of (i).}\, In the random graph, the expected
codegree is roughly $np^2 \gg 2$, but since we do not know how far $p$
exceeds $n^{-1/2}$, we need a slightly more careful argument.  Let $S$
be the sum of the codegrees $\sum_{\{u, v\}} d(u, v)$ over all
unordered pairs $\{u, v\}$, and let us decompose $S = S_1 + S_2 +
S_3$, where $S_1$ is the contribution from summands with $d(u, v) \in
\{0, 1\}$, $S_2$ is the contribution from summands with $2 \leq d(u,
v) \leq 4np^2$, and $S_3$ is the remainder.  We aim to show that $S_2
= \Omega(n^3 p^2)$, which will imply the result.

By double-counting, $S = \sum_v {d(v) \choose 2}$, where $d(v)$ is the
degree of vertex $v$.  By convexity, this is always at least $n {d
  \choose 2}$, where $d$ is the average degree.  Since $G_m$ has
exactly $m$ edges, $d = 2m/n = np \gg 1$.  Therefore, $S
\geq (0.5+o(1)) n(np)^2$.

On the other hand, Lemma \ref{lem:codegree-sequence} shows that \whp,
$G_m$ has the property that all codegrees are at most $np^2 \log n$,
and for every integer $4 \leq k \leq \log n$, the number of pairs with
codegree at least $knp^2$ is at most $n^2/k^3$.  Conditioning on
this, we may then bound $S_3$, the sum of codegrees which exceed
$4np^2$, by:
\begin{eqnarray*}
  S_3 &\leq& \sum_{k=4}^{\log n} (k+1)np^2 \cdot \frac{n^2}{k^3} \\
  &\leq& \frac{5}{4} \sum_{k=4}^{\log n}  \frac{n^3 p^2}{k^2} \\
  &\leq& n^3 p^2 \cdot \frac{5}{4}\left(\frac{\pi^2}{6} - \frac{1}{1^2} - \frac{1}{2^2} - \frac{1}{3^2}\right) \\
  &\leq& 0.4 n^3 p^2.
\end{eqnarray*}
Also, $S_1$, the sum of codegrees which are in $\{0, 1\}$, is
trivially at most ${n \choose 2} \ll n^3 p^2$ since we assumed $p \gg
n^{-1/2}$.  So, $S_2$, the sum of codegrees between 2 and
$4np^2$, is at least $ S_2 = S - S_1 - S_3 \geq 0.05 n^3 p^2$.
Therefore, \whp\ the number of pairs with codegree at least 2 is at least 
$0.05 n^3 p^2/(4np^2) = \Omega(n^2)$, as claimed. \hfill $\Box$

\vspace{2mm}

\noindent {\bf Proof of (ii).}\, The second part follows from Theorem
\ref{thm:balanced-extension}(ii) because $(\hme{K_4}{}, K_4)$ is a
balanced extension pair and $n^2 p^5$ is a negative power of $n$.
Let us now concentrate on the first part.  Conditioning on the high
probability event in claim (i), we may now assume that in $G_m$, the
proportion of pairs with codegree at least 2 is some $\gamma = \Omega(1)$.
Consider the $(i+1)$-st round, where $m \leq i < 2m$.  Regardless of
the choice of strategy, if all three incoming edges span pairs that
each have codegree at least 2, then we will be forced to
create a new copy of $\hme{K_4}{}$.  Incoming edges are uniformly
distributed over unoccupied pairs, and the number of occupied pairs in
$G_i$ is exactly $i = o(n^2)$.  So, since $G_i \supset G_m$, the
probability that all three incoming edges span pairs with codegree at
least 2 is $q \geq (1+o(1))\gamma^3 = \Omega(1)$.

Let $i$ run from $m$ to $2m$.  Then, the number of copies of
$\hme{K_4}{}$ in $G_{2m}$ is at least $\bin(m, q)$.  By the Chernoff bound,
this exceeds $mq/2 = \Omega(n^2 p)$ \wep, so we are
done. \hfill $\Box$

\vspace{2mm}

\noindent {\bf Proof of (iii).}\, Consider the $(i+1)$-st round, where
$2m \leq i < 3m$.  Regardless of the choice of strategy, if all three
incoming edges span pairs that complete copies of $K_4$, we will lose.
We can lower bound the number of such pairs by $\Omega(n^2 p)$ by
conditioning on the following events.  Let $A$ be the event that
$G_{2m}$ contains $\Omega(n^2 p)$ copies of $\hme{K_4}{}$, which occurs
\whp\ by (ii).  Also by (ii), with probability $1 - o(n^{-2})$,
every pair of vertices in $G_i \subset G_{3m}$ extends $O(1)$ copies of $\hme{K_4}{}$ into
$K_4$; call this event $B_i$.

Even after conditioning, the incoming edges in this round are still
independently and uniformly distributed over the ${n \choose 2} - i =
\Theta(n^2)$ unoccupied pairs of $G_i$.  Therefore, the probability
that both pairs complete $K_4$, conditioned on survival through the
$i$-th round, is $p_i = \Omega\big( \big(\frac{n^2 p}{n^2}\big)^3 \big) = \Omega(p^3)$.
Letting $i$ run from $2m$ to $3m$, we see that the probability that
any strategy can survive for $3m$ rounds is at most
\begin{eqnarray*}
  \mathbb{P} &\leq& \pr{\neg A} + \sum \pr{\neg B_i} + \prod (1-p_i) 
  \leq o(1)
  + \exp\left\{-\sum p_i\right\} \\
  &\leq& o(1) + \exp\Big\{ -\Omega(n^2 p \cdot p^3) \Big\} =
  o(1) + e^{-\omega(1)} \,= \, o(1),
\end{eqnarray*}
which completes the proof. \hfill $\Box$


\vspace{4mm}

It remains to establish Lemma \ref{lem:codegree-sequence}, which we
used to control the distribution of codegrees in claim (i) of the
upper bound.

\vspace{4mm}

\noindent \textbf{Proof of Lemma \ref{lem:codegree-sequence}.}\, Each
codegree is distributed as $\bin(n-2, p^2)$, so a union bound shows
that the probability that some codegree exceeds $np^2 \log n$ is at
most
\begin{displaymath}
  \mathbb{P} 
  \leq n^2 \cdot {n \choose np^2 \log n} (p^2)^{np^2 \log n} 
  \leq n^2 \cdot \left(\frac{enp^2}{np^2 \log n}\right)^{np^2 \log n}
  = o(1).
\end{displaymath}
Next, fix any $4 \leq k \leq \log n$, and let $X$ be the number of
pairs with codegree at least $knp^2$.  Consider an arbitrary vertex
$v$, and let $X_v$ be the number of vertices $u \neq v$ such that the
codegree of $\{v, u\}$ is at least $knp^2$.  Note that $X =
\frac{1}{2} \sum X_v$.

Since $d(v)$ is binomially distributed $\bin[n-1, p]$ and 
$np$ is a positive power of $n$, the degree $d(v)$ is at most
$1.1np$ \wep\ by Chernoff.  Condition on this, and condition further on a
neighborhood $N(v)$ of size $d(v)$.  For each $w \not \in N(v) \cup
\{v\}$, define the indicator random variable $I_w$ to be 1 if and only
if the codegree of $\{v, w\}$ is at least $knp^2$, or equivalently, if
$w$ has at least $knp^2$ neighbors in $N(v)$.  Note that because we
already fixed $N(v)$, these $I_w$ are independent since they are determined
by disjoint sets of edges.  Yet
$k \geq 4$ and $np^2 \gg 1$, so each $I_w$ has probability
\begin{displaymath}
  q = \pr{I_w}
  \leq {1.1np \choose knp^2} p^{knp^2} 
  \leq \left( \frac{1.1enp^2}{knp^2} \right)^{knp^2}
  \leq \left( \frac{3}{k} \right)^{knp^2}
  \ll \frac{1}{k^3}.
\end{displaymath}
Therefore, $X_v$ is stochastically dominated by $d(v) + \bin[n-1-d(v),
q]$.  Since $k \leq \log n$, a Chernoff bound implies that \wep, $X_v
\leq 1.1np + 2nq = o(n/k^3)$, which gives $X = \frac{1}{2} \sum X_v =
o(n^2/k^3)$.  The result follows by taking a union bound over all $v$
and $4 \leq k \leq \log n$.  \hfill $\Box$

\section{Concluding remarks}

\begin{itemize}
\item Although our theorems treat specific graphs (cycles, cliques,
  and complete bipartite graphs), we conjecture that the thresholds
  for avoiding general graphs $H$ follow from the natural
  generalization of the recipe that we used.

  To apply our machinery from Section \ref{sec:abstract-argument}, the
  first thing that we needed to specify was the parameter $s$.  This
  was the number of levels of danger considered by the avoidance
  strategy in the proof of the lower bound.  The correct choice of $s$
  then determined $\theta$, the negative exponent in the threshold (in
  terms of $p$) for avoidance:
  \begin{displaymath}
    \theta(H, r, s) = \frac{r^s (v(H)-2) + 2}{r^s\left(e(H) - s\right) + \frac{r^s - 1}{r-1}}.
  \end{displaymath}
  Furthermore, it is clear that the threshold for avoiding any fixed
  subgraph $H' \subset H$ is a lower bound for the threshold for
  avoiding $H$ itself.  This is because any strategy that avoids $H'$
  will certainly avoid $H$ as well.

  In light of this, we conjecture that the threshold for avoiding $H$
  in the Achlioptas process with parameter $r$ is $n^{2-\theta^*}$,
  where $\theta^*$ is the minimum value of $\theta(H', r, s)$ when $s$
  runs over all nonnegative integers and $H'$ runs over all subgraphs
  of $H$.

\item Just as in the case of analyzing the Achlioptas process for
giant component avoidance \cite{BKi}, one can also consider
the {\em offline} version of the fixed subgraph avoidance problem. In this
offline version, all random $r$-tuples of edges arriving during the
process are accessible to an algorithm, and it can make its choices
at each round, relying on the perfect knowledge of the past and the
future. The question is still how long the algorithm can typically avoid
the appearance of a copy of a fixed graph $H$. We expect that in
most of the cases there will be a sizable difference between the
online and the offline thresholds. Here is a sketch of the
illustrative case of $H=K_3$, $r=2$. For this case we can prove
that if $m=o(n^{4/3})$, then one can \whp\ avoid a copy of $K_3$
during the first $m$ rounds in the offline version. This should be
compared to the threshold of $m=n^{6/5}$ for the online version, given
by Theorem \ref{thm:main}. The argument proceeds as follows. Set
$p=2m/n^2$. The offline model in this case can be approximated quite
accurately by generating a random graph $G$ according to the
distribution $G(n,2m)$, and then splitting the edges of $G$
randomly into $m$ pairs: $(e_1,f_1),\ldots, (e_m,f_m)$. Denote the
above random matching of $E(G)$ by $\pi$.
We use the following strategy,
while processing the pairs $(e_i,f_i)$: in each pair $(e_i,f_i)$
choose an arbitrary edge not participating in any triangle in $G$,
otherwise pick an arbitrary edge. It is obvious that using
this strategy we can only lose (i.e. create a triangle) if $G$
contains a triangle with edges $x_1,x_2,x_3$ such that their
respective pairings in $\pi$ also belong to triangles in $G$.   
The number of triangles in $G$ is \whp\ of order $n^3p^3$, and
therefore the probability of having a triangle whose three edges
are paired in $\pi$ with edges from triangles is at most of order
$$
n^3p^3 \left(\frac{n^3p^3}{n^2p}\right)^3=n^6p^9=o(1).
$$
It would be very interesting to obtain tight results for the
offline small subgraph avoidance version of the Achlioptas process
for a wide variety of graphs $H$ and parameter $r$. 

\item The appearances of the giant component and of a fixed graph
are just two instances that have been addressed so far in the
context of the Achlioptas process. Naturally, one
can consider other graph theoretic properties as well in this
context. We hope to return to questions of this type in the
future.

\end{itemize}

\vspace{0.2cm}
\noindent {\bf Acknowledgment.}\, The authors would like to thank the
referee for suggestions that improved the exposition of this
paper.

\appendix

\section{Supporting results for avoiding $\boldsymbol{K_t}$}

In this section, we collect the supporting lemmas and inequalities
that are used to prove thresholds for avoiding $K_t$.
\shortversion{%
  Following a suggestion of the referee to shorten this paper, we do
  not provide complete proofs of all of these results.  Rather, we
  stop once each statement has been reduced to an inequality in
  several variables.  At that point, the remaining analysis is not so
  interesting, because such statements can of course can in theory be
  verified (although efficient proofs of non-polynomial inequalities
  in up to eight variables are not necessarily routine).  The interested
  reader can find the complete proofs on the arXiv at
  \texttt{http://arxiv.org/abs/0708.0443}.
}%
{%
  We will make claims about properties of simple single-variable
  functions without proof, as they can be verified by routine
  analytical methods.
}%

Throughout the appendix, we will set $s = \lfloor \log_r [(r-1)t + 1]
\rfloor$.  We begin by proving some basic facts about $s$.

\begin{lemma}
  \label{lem:s-dec-in-r}
  For fixed $t \geq 3$, the parameter $s$ is decreasing in $r$ in the
  range $r \geq 2$.
\end{lemma}

\noindent {\bf Proof.}\, %
\shortversion{%
  This follows by routine calculus, as it is not difficult to show
  that $\frac{\partial f}{\partial r} < 0$.  \hfill $\Box$
}{%
Let $t$ be fixed, and consider the function
$f(r) = \frac{\log [(r-1)t + 1]}{\log r}$.  It suffices to show that
$\frac{\partial f}{\partial r} < 0$.  Calculation yields:
\begin{displaymath}
  \frac{\partial f}{\partial r} 
  = \frac{(\log r)\frac{t}{(r-1)t+1} - \log[(r-1)t+1] \frac{1}{r}}{\log^2 r}
\end{displaymath}
We will show that the numerator is negative.  Since $t \geq
3$ and $r \geq 2$, $(r-1)t+1 \geq 3r-2 \geq 2r$, so the numerator
is less than
\begin{eqnarray*}
  \text{numerator} &<& (\log r)\frac{1}{r-1} - (\log 2r)\frac{1}{r} \\
  &=& (\log r)\left(\frac{1}{r-1} - \frac{1}{r}\right) - (\log 2) \frac{1}{r} \\
  &=& \left(\frac{\log r}{r-1} - \log 2\right) \frac{1}{r}.
\end{eqnarray*}
Now $\frac{\log r}{r-1}$ is a decreasing function of $r$ when $r \geq
2$, so it is at most $\log 2$; therefore, the entire expression is
always negative, and $f$ is indeed decreasing in $r$. \hfill $\Box$
}

\vspace{4mm}

\begin{lemma}
  \label{lem:s-leq-t/2}
  If $t \geq 4$ and $r \geq 2$, then $s \leq t/2$.  Furthermore, if $t
  \geq 5$ and $r \geq 2$, or if $t = 4$ and $r \geq 4$, then $s <
  t/2$.
\end{lemma}

\noindent {\bf Proof.}\, By Lemma \ref{lem:s-dec-in-r}, if $r \geq 2$,
then $s \leq \lfloor \log_2 (t+1) \rfloor$, and one may verify that
this is in turn $\leq t/2$ for all $t \geq 4$, and $< t/2$ for $t \geq
5$.  For the other range, when $r \geq 4$, Lemma \ref{lem:s-dec-in-r}
gives $s \leq \lfloor \log_4 (3t+1) \rfloor$, which is less than $t/2$
at $t=4$.  This finishes the lemma.  \hfill $\Box$

\subsection{Balanced graphs and extensions}

\begin{lemma}
  \label{lem:kt-lower-i-balanced-graph}
  For any $t \geq 4$ and $r \geq 2$,
  $\hme{K_t}{s}$ is a balanced graph.
\end{lemma}

\noindent {\bf Proof.}\, We must show that the edge density (number of
edges divided by number of vertices) of $\hme{K_t}{s}$ is at least as
large as the edge density of any of its proper induced subgraphs.  The edge
density of $\hme{K_t}{s}$ is exactly $\big[{t \choose 2} - s\big]/t$.
Lemma \ref{lem:s-leq-t/2} established that $s \leq t/2$, so the edge
density is at least $\big[{t \choose 2} - \frac{t}{2}\big]/t =
{t-1 \choose 2}/(t-1)$.  Yet the final quantity is precisely
the edge density of $K_{t-1}$, which is an upper bound on the edge
density of any proper induced subgraph of any $t$-vertex graph, so we are
done.  \hfill $\Box$

\vspace{4mm}

\begin{lemma}
  \label{lem:kt-lower-i-balanced-extension}
  For any $t \geq 4$ and $r \geq 2$, $\hme{K_t}{s}$ has the balanced
  extension property.
\end{lemma}

\noindent {\bf Proof.}\, Fix any graph $G$ of the form $\hme{K_t}{s}$,
and let $u$, $v$ be any two nonadjacent vertices of $G$.  We must show
that the function $e(H)/(v(H)-2)$ is maximal at $H = G$, where $H$ is
allowed to range over all induced subgraphs of $G$ that contain $\{u,
v\}$.  For any graph $H$ with $n$ vertices that is missing at least
one edge (e.g., the edge between $\{u, v\}$), $e(H)/(v(H)-2) \leq
\big[{n \choose 2} - 1\big]/(n-2) = (n+1)/2$.  For any proper induced
subgraph $H \subset G$, we then have $e(H)/(v(H)-2) \leq t/2$.  

Yet $e(G)/(v(G)-2) = \big[{t \choose 2} - s\big]/(t-2)$, and Lemma
\ref{lem:s-leq-t/2} established that $s \leq t/2$.  Using this bound
for $s$, we see that $e(G)/(v(G)-2) \geq \big[{t \choose 2} -
\frac{t}{2}\big]/(t-2) = t/2$, which matches our upper bound for
$e(H)/(v(H)-2)$, so we are done. \hfill $\Box$

\vspace{4mm}

\begin{lemma}
  \label{lem:kt-upper-balanced-extension}
  Suppose that $t \geq 4$.  Let $H_1 = K_{\lfloor \frac{t}{2} \rfloor,
    \lceil \frac{t}{2} \rceil}$, and arbitrarily choose the rest of the sequence $\{H_2, H_3,
  \ldots, H_f\}$, where $H_f = K_t$, by adding one edge at a time.
  Then every consecutive pair $(H_k, H_{k+1})$ is a balanced extension
  pair.
\end{lemma}

\noindent {\bf Proof.}\, Consider a consecutive pair $(H_k, H_{k+1})$.
By the construction, $H_k$ contains a complete bipartite subgraph that
was $H_1$; let $V_1 \cup V_2$ be the corresponding partition of the
vertex set.  Let $u$ and $v$ be the endpoints of the edge on which
$H_k$ and $H_{k+1}$ differ.  Without loss of generality, suppose that
$u, v \in V_1$.  (They must lie in the same part because $H_k$ already
contains all edges between $V_1$ and $V_2$.)  Now, consider any
subsets $U_1 \subset V_1$ and $U_2 \subset V_2$ such that $u, v \in
U_1$ and $U_1 \cup U_2 \neq V_1 \cup V_2$.  Let $H_k'$ be the subgraph
of $H_k$ induced by $U_1 \cup U_2$.  It suffices to show that
$e(H_k)/(v(H_k)-2) \geq e(H_k')/(v(H_k')-2)$.

Let us denote $u_1 = |U_1|$, $u_2 = |U_2|$, and let $e_1$ and $e_2$ be
the respective numbers of edges of $H_k$ spanned by $U_1$ and by
$U_2$.  Since the number of edges between $V_1$ and $V_2$ is
$\big\lfloor \frac{t}{2} \big\rfloor \big\lceil \frac{t}{2} \big\rceil
= \big\lfloor \frac{t^2}{4} \big\rfloor$, the number of edges in $H_k$
is at least $e_1 + e_2 + \big\lfloor \frac{t^2}{4} \big\rfloor$.  On
the other hand, the number of edges in $H_k'$ is precisely $e_1 + e_2
+ u_1 u_2$.  Thus, \shortversion{the result follows from the inequality below (proved in the full version):}{it suffices to show that:}
\begin{displaymath}
  \frac{e_1 + e_2 + \big\lfloor \frac{t^2}{4} \big\rfloor}{t-2} \geq \frac{e_1 + e_2 + u_1 u_2}{u_1 + u_2 - 2}.
\end{displaymath}
\shortversion{}{%
Since the quantity $(e_1 + e_2)$ appears in both numerators but the
denominators satisfy $t-2 > u_1 + u_2 - 2$, it will only strengthen
the inequality to increase $(e_1 + e_2)$ to its maximum possible
value, ${u_1 \choose 2} + {u_2 \choose 2} - 1$ (we must subtract 1
because we assumed that $H_k$ has no edge between $u, v \in U_1$).
Hence it suffices to show
\begin{displaymath}
  \frac{{u_1 \choose 2} + {u_2 \choose 2} - 1 + \big\lfloor \frac{t^2}{4} \big\rfloor}{t-2} \geq \frac{{u_1 \choose 2} + {u_2 \choose 2} - 1 + u_1 u_2}{u_1 + u_2 - 2}.
\end{displaymath}
Let $x = u_1 + u_2$.  The right hand side then simplifies to
$(x+1)/2$.  Also, ${u_1 \choose 2} + {u_2 \choose 2} \geq {x \choose
  2} - \big\lfloor \frac{x^2}{4} \big\rfloor = \big\lceil
\frac{x^2}{4} - \frac{x}{2} \big\rceil$.  Therefore, it suffices to
show the following inequality, which is equivalent to the next two
inequalities:
\begin{eqnarray*}
  \frac{\big\lceil \frac{x^2}{4} - \frac{x}{2} \big\rceil - 1 + \big\lfloor \frac{t^2}{4} \big\rfloor}{t-2} &\geq& \frac{x+1}{2} \\
  \Big\lceil \frac{x^2}{4} - \frac{x}{2} \Big\rceil - 1 + \Big\lfloor \frac{t^2}{4} \Big\rfloor &\geq& \frac{(x+1)(t-2)}{2} \\
  \Big\lceil \frac{x^2}{4} + \frac{x}{2} (1-t) - \frac{t}{2} + \Big\lfloor \frac{t^2}{4} \Big\rfloor \Big\rceil &\geq& 0.
\end{eqnarray*}
The expression under the ceiling sign is a quadratic polynomial in
$x$, so it is minimized at $x = t-1$, where it takes the value
$\big\lfloor \frac{t^2}{4} \big\rfloor - \frac{t^2}{4} - \frac{1}{4}$.
Since only 0 and 1 are quadratic residues modulo 4, this is always at
least $-1/2$, so we are done.} \hfill $\Box$

\subsection{Inequalities}

For the reader's convenience, we reproduce the definitions of the
parameters $s$ and $\theta$:
\begin{displaymath}
  s = \lfloor \log_r [(r-1)t + 1] \rfloor, \quad \quad \quad \theta = \frac{r^s (t-2) + 2}{r^s\left[{t \choose 2} - s\right] + \frac{r^s - 1}{r-1}}.
\end{displaymath}
\shortversion{The following inequalities are proved in the full
  version of this paper.}{}

\begin{inequality}
  \label{ineq:kt-upper-t/2-codegree}
  Suppose that either $t \geq 5$ and $r \geq 2$, or $t = 4$ and $r
  \geq 4$.  Then $-\theta \geq -\big\lfloor \frac{t}{2}
  \big\rfloor^{-1}$.
\end{inequality}

\shortversion{}{\noindent {\bf Proof.}\, It suffices to prove the stronger inequality
that $\theta \cdot \frac{t}{2} \leq 1$.  Since the denominator of $\theta$
is positive, this is equivalent to
\begin{displaymath}
  \frac{t}{2} r^s (t-2) + t \leq r^s \left[\frac{t}{2} (t-1) - s\right] + \frac{r^s - 1}{r-1}.
\end{displaymath}
Rearranging terms, this is equivalent to
\begin{equation}
  \label{ineq:i:kt-u-t-c:0}
  t \leq r^s \left(\frac{t}{2} - s\right) + \frac{r^s - 1}{r-1}.
\end{equation}
Now by the definition of $s = \lfloor \log_r [(r-1)t + 1] \rfloor$, we
must have $r^{s+1} > (r-1)t + 1$.  Since $r$, $s$, and $t$ are integers,
this implies $r^{s+1} \geq (r-1)t + 2$, from which we deduce (using $r \geq 2$):
\begin{equation}
  \label{ineq:i:kt-u-t-c:11}
  r^s \geq \frac{(r-1)t + 2}{r} = t - \frac{t-2}{r} \geq t - \frac{t-2}{2} = \frac{t+2}{2}.
\end{equation}
Since $r$ and $s$ are integers, we in fact have
\begin{equation}
  \label{ineq:i:kt-u-t-c:1}
  r^s \geq \left\lceil \frac{t+2}{2} \right\rceil.
\end{equation}
Also, since $s$ is always $\geq 1$ for $t \geq 1$,
\begin{equation}
  \label{ineq:i:kt-u-t-c:2}
  \frac{r^s - 1}{r-1} \geq 1.
\end{equation}

Let us now split into several cases, depending on the (integral) values of $t$ and $r$.
\begin{description}
\item[Case 1: $t \geq 10$.] Let us bound the terms on the right hand
  side of inequality \eqref{ineq:i:kt-u-t-c:0} one at a time.  By
  \eqref{ineq:i:kt-u-t-c:1}, $r^s > t/2$, and by
  \eqref{ineq:i:kt-u-t-c:2}, $\frac{r^s - 1}{r-1} > 0$.  Also, Lemma
  \ref{lem:s-dec-in-r} showed that $s$ was decreasing in $r$, so for
  $r \geq 2$, we have $s \leq \lfloor \log_2 (t+1) \rfloor$.  One may
  routinely verify that $\lfloor \log_2 (t+1) \rfloor \leq \frac{t}{2}
  - 2$ for $t \geq 10$.  Combining all of these bounds, we obtain
  \begin{displaymath}
    r^s \left(\frac{t}{2} - s\right) + \frac{r^s - 1}{r-1} > \left(\frac{t}{2}\right) (2) + 0 = t,
  \end{displaymath}
  as desired.

\item[Case 2: $t \geq 5$ and $r \geq 6$.] We follow the same method as
  the previous case.  If we use $r \geq 6$ in the derivation of
  inequality \eqref{ineq:i:kt-u-t-c:11}, we obtain
  \begin{displaymath}
    r^s \geq t - \frac{t-2}{r} \geq t - \frac{t-2}{6} = \frac{5t+2}{6} > \frac{5}{6} t.
  \end{displaymath}
  By \eqref{ineq:i:kt-u-t-c:2}, $\frac{r^s - 1}{r-1} > 0$.  Also, by
  Lemma \ref{lem:s-dec-in-r}, if $r \geq 6$, then $s \leq \lfloor
  \log_6 (5t+1) \rfloor$.  One may routinely verify that $\lfloor
  \log_6 (5t+1) \rfloor \leq \frac{t}{2} - \frac{6}{5}$ for $t \geq
  5$.  Combining all of these bounds, we obtain
  \begin{displaymath}
    r^s \left(\frac{t}{2} - s\right) + \frac{r^s - 1}{r-1} > \left(\frac{5}{6} t\right) \left(\frac{6}{5}\right) + 0 = t,
  \end{displaymath}
  as desired.

\item[Case 3: $t = 4$ and $r \geq 4$.] By Lemma
  \ref{lem:s-dec-in-r}, if $r \geq 4$, then $s \leq \lfloor \log_4
  (3t+1) \rfloor$.  Since we have $t = 4$, this means that $s \leq 1$,
  and in particular, $\frac{t}{2} - s \geq 1$.  Yet inequality
  \eqref{ineq:i:kt-u-t-c:1} bounds $r^s \geq \lceil (t+2)/2 \rceil =
  3$ and inequality \eqref{ineq:i:kt-u-t-c:2} bounds $\frac{r^s - 1}{r-1}
  \geq 1$, so combining all of these bounds, we obtain
  \begin{displaymath}
    r^s \left(\frac{t}{2} - s\right) + \frac{r^s - 1}{r-1} \geq 3 (1) + 1 = 4 = t,
  \end{displaymath}
  as desired.

\item[Case 4: remainder.] Only finitely many instances remain: when
  the integers $t$ and $r$ satisfy $5 \leq t \leq 9$ and $2 \leq r
  \leq 5$.  A routine check verifies that $\theta \cdot \frac{t}{2} \leq 1$
  in all 20 of these cases.

\end{description}

This completes the proof.  \hfill $\Box$

\vspace{4mm}

}  

\begin{inequality}
  \label{ineq:kt-upper-pos-power}
  For any $t \geq 4$ and $r \geq 2$, if $p \gg n^{-\theta}$, then
  $n^{t-2} p^{{t \choose 2} - s - 1}$ is a positive power of $n$.
\end{inequality}

\shortversion{}{\noindent {\bf Proof.}\, Clearing the denominator of $\theta$, we see that
we must establish the following inequality:
\begin{displaymath}
  (t-2)\left[ r^s\left[{t \choose 2} - s\right] + \frac{r^s - 1}{r-1} \right] > \left[{t \choose 2} - s - 1\right] \left[r^s (t-2) + 2\right].
\end{displaymath}%
The right hand side is equal to 
\begin{displaymath}
  \left[{t \choose 2} - s - 1\right] \left[r^s (t-2) + 2\right] = \left[{t \choose 2} - s\right] r^s(t-2) + \left[{t \choose 2} - s\right] 2 - r^s(t-2) - 2,
\end{displaymath}
so our initial inequality is equivalent to
\begin{displaymath}
  (t-2) \frac{r^s - 1}{r-1} > \left[{t \choose 2} - s\right] 2 - r^s (t-2) - 2,
\end{displaymath}
so by rearranging terms and ignoring the $s$ term appearing in the
bracket on the right hand side, it suffices to establish
\begin{displaymath}
  (t-2)\frac{r^{s+1} - 1}{r-1} \geq 2{t \choose 2} - 2.
\end{displaymath}
Dividing through by $t-2$, we see that it suffices to show that
\begin{displaymath}
  \frac{r^{s+1} - 1}{r-1} \geq t+1.
\end{displaymath}
Since both sides of this inequality are integers, it is equivalent to
\begin{displaymath}
  \frac{r^{s+1} - 1}{r-1} > t.
\end{displaymath}
Yet by definition of $s = \lfloor \log_r [(r-1)t + 1] \rfloor$, we
must have $r^{s+1} > (r-1)t + 1 \Rightarrow \frac{r^{s+1} - 1}{r-1} >
t$, which completes the proof. \hfill $\Box$

\vspace{4mm}

} 

\begin{inequality}
  \label{ineq:kt-upper-neg-power}
  Suppose that $t \geq 5$ and $r \geq 2$, or $t = 4$ and $r \geq 4$.
  If $p \ll n^{-\theta}$, then $n^{t-2} p^{{t \choose 2} - s}$ is a
  negative power of $n$.
\end{inequality}

\shortversion{}{\noindent {\bf Proof.}\, Clearing the denominator of $\theta$, we see
that we must establish the following inequality:
\begin{displaymath}
  (t-2)\left[ r^s\left[{t \choose 2} - s\right] + \frac{r^s - 1}{r-1} \right] < \left[{t \choose 2} - s\right] \left[r^s (t-2) + 2\right].
\end{displaymath}
This is equivalent to
\begin{displaymath}
  (t-2)\frac{r^s - 1}{r-1} < \left[{t \choose 2} - s\right] 2.
\end{displaymath}
By definition of $s = \lfloor \log_r [(r-1)t + 1] \rfloor$, we must
have $r^s \leq (r-1)t+1$.  Therefore, the left hand side is at most
$(t-2)t$, so it suffices to show that
\begin{displaymath}
  t^2 - 2t < t^2 - t - 2s.
\end{displaymath}
But this final inequality is a consequence of Lemma
\ref{lem:s-leq-t/2}, which established that $s < t/2$. \hfill $\Box$

}  

\section{Supporting results for avoiding $\boldsymbol{K_{t,t}}$}

Coincidentally, the definition of the parameter $s$ is exactly the same
for avoiding $K_t$ and avoiding $K_{t,t}$, so we can still use Lemmas
\ref{lem:s-dec-in-r} and \ref{lem:s-leq-t/2} (which prove properties
of $s$) in this section.  The specification of $\theta$ will be
different, however.  For the reader's convenience, we reproduce the
definitions here.
\begin{displaymath}
  s = \lfloor \log_r [(r-1)t + 1] \rfloor, \quad\quad\quad \theta = \frac{r^s (2t-2) + 2}{r^s (t^2 - s) + \frac{r^s - 1}{r-1}}.
\end{displaymath}

\subsection{Balanced graphs}

\begin{lemma}
  \label{lem:ktt-lower-i-balanced-graph}
  For any $t \geq 3$ and $r \geq 2$, $\hme{K_{t,t}}{s}$ is a balanced
  graph, except in the case when $t=3$, $r=2$, and the graph is
  $K_{2,3}$ with a pendant edge.  In that final case, if $p \gg
  n^{-18/31} / \log n$, the number of copies of that graph in $G_m$ is
  still $O(n^6 p^7)$ \wep.
\end{lemma}

\noindent {\bf Proof.}\, We must show that the edge density (number of
edges divided by number of vertices) of $\hme{K_{t,t}}{s}$ is at least
the edge density of any proper induced subgraph.  The edge density of the
complete bipartite graph $K_{a,b}$ is $ab/(a+b)$, which is increasing
in both $a$ and $b$, so the edge density of any proper induced subgraph of
$\hme{K_{t,t}}{s}$ is at most $t(t-1)/(2t-1)$.  On the other hand, the
edge density of $\hme{K_{t,t}}{s}$ is precisely $(t^2 - s)/(2t)$, so
we must show that
\begin{displaymath}
  \frac{t(t-1)}{2t-1} \leq \frac{t^2 - s}{2t}.
\end{displaymath}
Clearing the denominators, this is equivalent to
\begin{displaymath}
  2t^3 - 2t^2 \leq 2t^3 - t^2 - s(2t-1).
\end{displaymath}
Rearranging terms, this is equivalent to
\begin{displaymath}
  s \leq \frac{t^2}{2t-1}.
\end{displaymath}
Now if $t \geq 4$, Lemma \ref{lem:s-leq-t/2} bounds $s \leq t/2$,
which finishes the inequality.  

The only remaining case is $t=3$.  However, Lemma \ref{lem:s-dec-in-r}
established that the dependence of $s = \lfloor \log_r [(r-1)t + 1]
\rfloor$ on $r$ was decreasing, so $s = 1$ for $r \geq 3$, and $s = 2$
for $r = 2$.  One may manually verify that of all of the graphs of the
form $\hme{K_{3,3}}{}$ and $\hme{K_{3,3}}{2}$, the only one which is
not balanced is the deletion from $K_{3,3}$ of two edges incident to
the same vertex, which is $K_{2,3}$ with a pendant edge, as claimed.
Since that graph, which we denote $K_{2,3} + e$, arises only when
$s=2$, this happens only when $r=2$.

Now let us bound the number of copies of that graph in $G(n, p)$, when
$p \gg n^{-\theta} / \log n$.  In the case $t=3, r=2$, we have
$\theta = -\frac{18}{31}$, and so $n^5 p^6$, roughly the expected number of
copies of $K_{2,3}$ in the random graph, is a positive power of $n$.
So, since $K_{2,3}$ is balanced, Theorem \ref{thm:balanced-graph}
bounds the number of copies of $K_{2,3}$ in $G_m$ by $O(n^5 p^6)$
\wep.  Also, $np$ is a positive power of $n$, so we may bound all
degrees by $2np$ \wep.  If both situations hold, we may conclude
that the number of copies of $K_{2,3} + e$ is $O(n^5 p^6 \cdot np) =
O(n^6 p^7)$, as desired.  \hfill $\Box$

\vspace{4mm}

\begin{lemma}
  \label{lem:ktt-lower-i-balanced-extension}
  For any $t \geq 3$ and $r \geq 2$, $\hme{K_{t,t}}{s}$ has the
  balanced extension property.
\end{lemma}

\noindent {\bf Proof.}\, Fix any graph $G$ of the form
$\hme{K_{t,t}}{s}$, and let $u$, $v$ be any two nonadjacent vertices
of $G$.  We must show that the function $e(H)/(v(H)-2)$ is maximized
at $H = G$, where $H$ is allowed to range over all proper induced subgraphs of $G$
that contain $\{u, v\}$.  Any such $H$ is still bipartite with respect
to $G$'s bipartition; suppose that it has $a$ vertices on one side and
$b$ on the other.  Since we assumed that $H$ is missing at least the
edge joining $\{u, v\}$, we must have $e(H)/(v(H)-2) \leq
(ab-1)/(a+b-2)$.  This is increasing in both $a$ and $b$, so its
maximum over proper induced subgraphs $H$ is $[t(t-1)-1]/(2t-3)$.  
\shortversion{Thus, the result follows from the inequality below
  (proved in the full version):}{So we must prove that}
\begin{displaymath}
  \frac{t^2 - s}{2t-2} \geq \frac{t(t-1)-1}{2t-3}.
\end{displaymath}
\shortversion{\hfill $\Box$}{%
Clearing the denominators and expanding brackets, this is equivalent to
\begin{displaymath}
  2t^3 - 3t^2 - s(2t-3) \geq 2t^3 - 4t^2 + 2.
\end{displaymath}
Rearranging terms, this is equivalent to
\begin{displaymath}
  t^2 - 2 \geq s(2t-3).
\end{displaymath}
If $t \geq 4$, then Lemma \ref{lem:s-leq-t/2} implies that $s \leq
t/2$, which implies the above inequality.  On the other hand, if
$t=3$, then Lemma \ref{lem:s-dec-in-r} established that $s \leq
\lfloor \log_2 (t+1) \rfloor$ for $r \geq 2$, so we must have $s \leq
2$, and this also implies the above inequality.  Therefore, we are
done.  \hfill $\Box$

} 

\vspace{4mm}

\begin{lemma}
  \label{lem:ktt-upper-i-balanced-extension}
  For any fixed positive integer $k$, consider the following 4-partite
  graph, which we call $H_1$.  Let the parts be $V_1, V_2, V_3, V_4$,
  with $|V_1| = |V_2| = k$ and $|V_3| = |V_4| = 1$, and place edges
  such that $(V_1, V_2)$, $(V_1, V_4)$, and $(V_3, V_2)$ are complete
  bipartite.  There are no more edges.  Let $H_2$ be obtained from
  $H_1$ by adding the edge between $V_3$ and $V_4$.
  Then $(H_1, H_2)$ is a balanced extension pair.
\end{lemma}

\noindent {\bf Proof.}\, Consider any subsets $U_1 \subset V_1$ and
$U_2 \subset V_2$, and let $H_1'$ be the subgraph of $H_1$ induced by
$U_1 \cup U_2 \cup V_3 \cup V_4$.  We must show that
$e(H_1')/(v(H_1')-2) \leq e(H_1)/(v(H_1)-2)$.  Let $a = |U_1|$ and $b
= |U_2|$.  Then, $\frac{e(H_1')}{v(H_1')-2} = \frac{ab+a+b}{a+b} =
\frac{ab}{a+b} + 1$, which is increasing in both $a$ and $b$.
Therefore, $\frac{e(H_1')}{v(H_1')-2} \leq \frac{k^2 + k + k}{k + k} =
\frac{e(H_1)}{v(H_1)-2}$, and we are done. \hfill $\Box$

\vspace{4mm}

\begin{lemma}
  \label{lem:ktt-upper-balanced-extension-even}
  Suppose that $t$ is even and at least $4$.  Let $H_1$ be the
  4-partite graph with parts $V_1, V_2, V_3, V_4$, each of size $t/2$,
  and edges such that $(V_1, V_2)$, $(V_1, V_4)$, and $(V_3, V_2)$ are
  complete bipartite.  Let $\{H_2, \ldots, H_{1+(t/2)}\}$ be obtained
  by successively adding single edges until $H_{1+(t/2)}$ has a
  perfect matching between $V_3$ and $V_4$.  Then, arbitrarily choose
  the rest of the sequence $\{H_{2+(t/2)}, \ldots, H_f\}$ by adding
  one edge at a time, until the final term is the complete bipartite
  graph $K_{t,t}$ with bipartition $(V_1 \cup V_3, V_2 \cup V_4)$.
  Then every consecutive pair $(H_k, H_{k+1})$ is a balanced extension
  pair.
\end{lemma}

The proof breaks into two cases, since there are two stages of edge
addition.  \shortversion{To give a flavor of the argument, we show how
  to reduce one of the cases to an inequality in several variables.}{}

\vspace{4mm}

\noindent {\bf Proof of Lemma \ref{lem:ktt-upper-balanced-extension-even}
  for $\boldsymbol{k \leq t/2}$.}\, Consider a consecutive pair $(H_k, H_{k+1})$.
By the construction, $H_k$ has the following structure.  The vertex
set is partitioned into $V_1 \cup V_2 \cup V_3 \cup V_4$, with all
parts of size $t/2$.  The pairs $(V_1, V_2)$, $(V_1, V_4)$, and $(V_3,
V_2)$ are complete bipartite graphs, and there is a $(k-1)$-edge
matching between $V_3$ and $V_4$.  There are no other edges.  Also,
there is a pair of vertices $u \in V_3$, $v \in V_4$, not involved in
the $(k-1)$-edge matching, at which the addition of an edge creates
$H_{k+1}$.  Now consider any family of subsets $U_i \subset V_i$ such
that $u \in U_3$ and $v \in U_4$.  Let $H_k'$ be the subgraph of $H_k$
induced by $\cup U_i$.  We must show that $e(H_k')/(v(H_k')-2) \leq
e(H_k)/(v(H_k)-2)$.

For brevity, let $a = |U_1|$, $b = |U_2|$, $c = |U_3|$, and $d =
|U_4|$.  Since the edges between $U_3$ and $U_4$ form a matching of at
most $k-1$ edges which does not involve $u \in U_3$ or $v \in U_4$,
there can be at most $\min\{c-1, d-1, k-1\} = \min\{c, d, k\} - 1$
edges there.  Therefore,
\begin{displaymath}
  \frac{e(H_k')}{v(H_k')-2} \leq \frac{ab + ad + cb + (\min\{c, d, k\} - 1)}{a+b+c+d-2}.
\end{displaymath}
\shortversion{%
  The result follows by showing that the right hand side is at most
  $\frac{\frac{3}{4}t^2 + (k-1)}{2t-2} = \frac{e(H_k)}{v(H_k)-2}$,
  which is done in the full version of this paper.  \hfill $\Box$
}%
{%
Let us simplify the bound by removing the variables $a$ and $b$ from
the inequality.  To do this, we show that the bound is nondecreasing in
both $a$ and $b$; then, we may replace them with their maximum value
$t/2$.  Calculating the partial derivative of the right hand side:
\begin{eqnarray*}
  \frac{\partial}{\partial a} &=& \frac{(a+b+c+d-2)(b+d) - (ab+ad+cb+ (\min\{c, d, k\} - 1))}{(a+b+c+d-2)^2} \\
  \text{numerator} &=& (b+d - 1)^2 + cd - \min\{c, d, k\} \\
  &\geq& cd - c,
\end{eqnarray*}
which is always nonnegative because $v \in U_4 \Rightarrow d \geq 1$.
An analogous result holds for $\frac{\partial}{\partial b}$ by
symmetry.  Hence we may substitute $a = b = t/2$ to obtain the bound
\begin{equation}
  \label{ineq:ktt-u-b-e:0}
  \frac{e(H_k')}{v(H_k')-2} 
  \leq \frac{ \frac{t^2}{4} + \frac{t}{2}(c+d) + \min\{c, d, k\} - 1 }{ t + (c+d) - 2 }
  = \frac{t}{2} + \frac{ -\frac{t^2}{4} + t + \min\{c, d, k\} - 1 }{ t + (c+d) - 2 }.
\end{equation}
We reserve the case $k = t/2$ for separate treatment at the end.
\begin{description}
\item[Case 1: $k < t/2$.] Using the bound $\min\{c, d, k\} \leq k$ in
  \eqref{ineq:ktt-u-b-e:0}, we obtain:
  \begin{displaymath}
    \frac{e(H_k')}{v(H_k')-2} 
    \leq \frac{t}{2} + \frac{ -\frac{t^2}{4} + t + k - 1 }{ t + (c+d) - 2 } 
  \end{displaymath}
  Note that if we substitute $c = d = t/2$ (their maximum values) in
  this bound, we obtain precisely $e(H_k)/(v(H_k)-2)$.  Therefore, it
  suffices to show that the bound is nondecreasing in $c$ and $d$.  We
  accomplish this by showing that the numerator of the large fraction
  is always $\leq 0$.

  We assumed in this case that $k < t/2$, but in fact, since $k$ is an
  integer and $t$ is even, we actually have $k \leq \frac{t}{2} - 1$.
  Substituting this into the numerator:
  \begin{displaymath}
    -\frac{t^2}{4} + t + k - 1 \leq -\frac{t^2}{4} + t + \left(\frac{t}{2} - 1\right) - 1 = -\frac{1}{4} [ (t-3)^2 - 1 ].
  \end{displaymath}
  Since we assumed that $t \geq 4$, this is always $\leq 0$, so the
  bound is indeed nondecreasing in $c$ and $d$.  This finishes this
  case.

\item[Case 2: $k = t/2$.]  This time, we bound $\min\{c, d, k\} \leq
  (c+d)/2$ in \eqref{ineq:ktt-u-b-e:0}.  Letting $x = (c+d)/2$, we
  have
  \begin{displaymath}
    \frac{e(H_k')}{v(H_k')-2} \leq \frac{t}{2} + \frac{ -\frac{t^2}{4} + t + x - 1 }{ t + 2x - 2 }.
  \end{displaymath}
  Note that if we substitute $x = t/2$ (its maximum value) in the
  bound, we obtain precisely $e(H_k)/(v(H_k)-2)$.  Therefore, it
  suffices to show that the bound is increasing in $x$.  Taking the
  numerator of the partial derivative of the fraction involving $x$:
  \begin{displaymath}
    \text{numerator of}\ \frac{\partial}{\partial x} = (t + 2x - 2) - \left( -\frac{t^2}{4} + t + x - 1 \right) 2 = \frac{t}{2} (t-2).
  \end{displaymath}
  Since we assumed $t \geq 4$, this is always positive, so the bound
  is indeed increasing in $x$.  This completes the proof. \hfill
  $\Box$
\end{description}

\vspace{4mm}

\noindent {\bf Proof of Lemma
  \ref{lem:ktt-upper-balanced-extension-even} for $\boldsymbol{k > t/2}$.}\,
Consider a consecutive pair $(H_k, H_{k+1})$.  By construction, $H_k$
has the following structure.  The vertex set is partitioned into $V_1
\cup V_2 \cup V_3 \cup V_4$, where each part has size $t/2$.  The
pairs $(V_1, V_2)$, $(V_1, V_4)$, and $(V_3, V_2)$ are complete
bipartite graphs, and there is a perfect matching between $V_3$ and
$V_4$.  There may be some more edges as well between $V_3$ and $V_4$,
but there are two vertices $u \in V_3$ and $v \in V_4$ such that there
is no edge between $u$ and $v$.  There are no more edges in the entire
graph.  Also, $H_{k+1}$ is obtained by $H_k$ by the addition of just
the edge joining $u$ and $v$.  Now, consider any family of subsets
$U_i \subset V_i$ such that $u \in U_3$ and $v \in U_4$.  Let $H_k'$
be the subgraph of $H_k$ induced by $\cup U_i$.  We must show that
$e(H_k')/(v(H_k')-2) \leq e(H_k)/(v(H_k)-2)$.

For brevity, let $a = |U_1|$, $b = |U_2|$, $c = |U_3|$, and $d =
|U_4|$.  Let $E$ be the number of edges between $U_3$ and $U_4$ in
$H_k$.  Then
\begin{displaymath}
  \frac{e(H_k')}{v(H_k')-2} = \frac{ab + ad + cb + E}{a+b+c+d-2}.
\end{displaymath}
Next, recall that $H_k$ contained a perfect $(t/2)$-edge matching
between $V_3$ and $V_4$.  The maximum number of edges of this matching
that are included in $E$ (i.e., go between $U_3$ and $U_4$) is
$\min\{c, d\} \leq (c+d)/2$.  Therefore, the number of edges in $H_k$ between $V_3$
and $V_4$ is at least $E + \frac{t}{2} - \frac{c+d}{2}$.  Since the
rest of $H_k$ consists of 3 complete bipartite subgraphs $K_{t/2, t/2}$:
\begin{displaymath}
  \frac{e(H_k)}{v(H_k)-2} \geq \frac{3 \frac{t^2}{4} + E + \frac{t}{2} - \frac{c+d}{2} }{2t-2},
\end{displaymath}
and so it suffices to show that
\begin{displaymath}
  \frac{ab + ad + cb + E}{a+b+c+d-2} \leq \frac{3 \frac{t^2}{4} + E + \frac{t}{2} - \frac{c+d}{2} }{2t-2}.
\end{displaymath}
Observe that $E$ appears in the numerator on both sides.  Since the
denominator on the left hand side is always less than or equal to the
denominator on the right hand side, we may replace $E$ (on both sides) by
its maximum possible value and only sharpen the inequality.  This
maximum is $|U_3| |U_4|-1 = cd-1$, because we assumed that $H_k$ had
no edge between $u \in U_3$ and $v \in U_4$.  Thus it suffices to show
that
\begin{displaymath}
  \frac{ab + ad + cb + (cd-1)}{a+b+c+d-2} \leq \frac{3 \frac{t^2}{4} + (cd-1) + \frac{t}{2} - \frac{c+d}{2} }{2t-2}.
\end{displaymath}
Note that the variables $a$ and $b$ appear only on the left hand side.  We
will show that the left hand side is nondecreasing in each of $a$ and $b$,
which will then allow us to replace each of them by their maximum
value ($t/2$) and only sharpen the inequality.  Let us verify this by
taking the partial derivative of the left hand side with respect to $a$:
\begin{eqnarray*}
  \frac{ab + ad + cb + (cd-1)}{a+b+c+d-2} &=& \frac{(a+c)(b+d)-1}{(a+c)+(b+d)-2} \\
  \frac{\partial}{\partial a} &=& \frac{[(a+c)+(b+d)-2](b+d) - [(a+c)(b+d) - 1]}{(a+b+c+d-2)^2} \\
  \text{numerator} &=& (b+d)^2 - 2(b+d) + 1,
\end{eqnarray*}
which is a perfect square, hence always nonnegative.  Thus the left hand
side was indeed nondecreasing in $a$.  By symmetry, the same is true
for $b$, and so it now suffices to prove:
\begin{displaymath}
  \frac{\frac{t^2}{4} + \frac{t}{2}(c+d) + (cd-1)}{t+(c+d)-2} \leq \frac{3 \frac{t^2}{4} + (cd-1) + \frac{t}{2} - \frac{c+d}{2} }{2t-2}.
\end{displaymath}
Now let $x = (c+d)/2$.  Then, $cd \leq x^2$.  Observe that $cd$
appears in both the numerator of the left hand side and the right hand side, but
the denominator on the left is at most the denominator on the right.
Therefore, we may increase $cd$ to $x^2$ on both sides, and only
sharpen the inequality.  Thus, it suffices to show that
\begin{displaymath}
  \frac{\frac{t^2}{4} + \frac{t}{2}(2x) + (x^2-1)}{t+(2x)-2} \leq \frac{3 \frac{t^2}{4} + (x^2-1) + \frac{t}{2} - x}{2t-2}.
\end{displaymath}
The left hand side simplifies:
\begin{displaymath}
  \frac{\frac{t^2}{4} + \frac{t}{2}(2x) + (x^2-1)}{t+(2x)-2} = \frac{\left(\frac{t}{2}+x\right)^2 - 1}{2\left[\left(\frac{t}{2} + x\right) - 1\right]} = \frac{\left(\frac{t}{2} + x\right) + 1}{2},
\end{displaymath}
so it suffices to show that
\begin{displaymath}
  \frac{\left(\frac{t}{2} + x\right) + 1}{2} \leq \frac{3 \frac{t^2}{4} + (x^2-1) + \frac{t}{2} - x}{2t-2}.
\end{displaymath}
Clearing the denominators by multiplying by $2t-2$, we see that it suffices to show that
\begin{displaymath}
  \left( \frac{t}{2} + x + 1 \right)(t-1) \leq 3\frac{t^2}{4} + (x^2-1) + \frac{t}{2} - x.
\end{displaymath}
Expanding the brackets and collecting all remaining terms on one side,
this is equivalent to
\begin{displaymath}
  0 \leq x^2 - xt + \frac{t^2}{4}.
\end{displaymath}
Recognizing the right hand side as the perfect square $(x -
\frac{t}{2})^2$, we are done. \hfill $\Box$

} 

\vspace{4mm}

\begin{lemma}
  \label{lem:ktt-upper-balanced-extension-odd}
  Suppose that $t$ is odd and at least 3.  Let $H_1$ be a 6-partite
  graph with parts $\{V_i\}_1^6$ such that $V_3$ and $V_4$ are
  singletons, and the other four parts each have size $\lfloor t/2
  \rfloor$.  Let there be edges be such that the two pairs $(V_1,
  V_2)$ and $(V_5, V_6)$ are each complete bipartite graphs, let the
  vertex in $V_3$ be adjacent to all of $V_2 \cup V_4 \cup V_6$, and
  let the vertex in $V_4$ be adjacent to all of $V_1 \cup V_3 \cup
  V_5$.  There are no more edges.

  Let $\{H_2, \ldots, H_{1 + \lfloor t/2 \rfloor}\}$ be obtained by
  successively adding single edges until $H_{1 + \lfloor t/2 \rfloor}$
  has a perfect matching between $V_1$ and $V_6$.  To create the next
  $\lfloor t/2 \rfloor$ graphs in the sequence, we put down a matching
  between $V_5$ and $V_2$, one edge at a time.  Finally, arbitrarily
  choose the rest of the sequence $\{H_{2 + 2 \lfloor t/2 \rfloor},
  \ldots, H_f\}$ by adding one edge at a time, until the final term is
  the complete bipartite graph $K_{t,t}$ with bipartition $(V_1 \cup
  V_3 \cup V_5, V_2 \cup V_4 \cup V_6)$.

  Then every consecutive pair $(H_k, H_{k+1})$ is a balanced extension
  pair.
\end{lemma}

The proof breaks into three cases, since there are three stages of
edge addition.  \shortversion{To give a flavor of the argument, we
  show how to reduce one of the cases to an inequality in several
  variables.}{}

\shortversion{}{

\vspace{4mm}

\noindent {\bf Proof of Lemma
  \ref{lem:ktt-upper-balanced-extension-odd} for $\boldsymbol{k \leq \lfloor t/2
  \rfloor}$.}\, Consider a consecutive pair $(H_k, H_{k+1})$.  By the
construction, $H_k$ has the following structure.  The vertex set is
partitioned into $\{V_i\}_1^6$, with $|V_3| = |V_4| = 1$ and all
other $|V_i| = \lfloor t/2 \rfloor$.  The pairs $(V_1, V_2)$ and
$(V_5, V_6)$ are each complete bipartite graphs, 
the vertex in $V_3$ is adjacent to all of $V_2 \cup V_4 \cup V_6$, 
the vertex in $V_4$ is adjacent to all of $V_1 \cup V_3 \cup V_5$, 
and there is a
$(k-1)$-edge matching between $V_1$ and $V_6$.  There are no other
edges.  Also, there is a pair of vertices $u \in V_1$, $v \in V_6$,
not involved in the $(k-1)$-edge matching, at which the addition of an
edge creates $H_{k+1}$.  Now consider any family of subsets $U_i
\subset V_i$ such that $u \in U_1$ and $v \in U_6$.  Let $H_k'$ be the
subgraph of $H_k$ induced by $\cup U_i$.  We must show that
$e(H_k')/(v(H_k')-2) \leq e(H_k)/(v(H_k)-2)$.

For brevity, let $a = |U_1|$, $b = |U_2|$, $c = |U_3|$, $d = |U_4|$,
$e = |U_5|$, and $f = |U_6|$.  Since the edges between $U_1$ and $U_6$
form a matching of at most $k-1$ edges which does not involve $u \in
U_1$ or $v \in U_6$, there can be at most $\min\{a-1, f-1, k-1\} =
\min\{a, f, k\} - 1$ edges there.  Therefore,
\begin{displaymath}
  \frac{e(H_k')}{v(H_k')-2} \leq \frac{ab + ef + c(b+f) + (a+e)d + cd + (\min\{a, f, k\} - 1)}{a+b+c+d+e+f-2}.
\end{displaymath}
We will show that the bound is nondecreasing in $c$, which will allow
us to replace $c$ with its maximum value of 1.  To do this, we as
usual consider the numerator of the partial derivative with respect to
$c$:
\begin{eqnarray*}
  \text{numerator of } \frac{\partial}{\partial c} &=& (a+b+c+d+e+f-2)(b+f+d) \\
  && \quad -\ [ab + ef + c(b+f) + (a+e)d + cd + (\min\{a, f, k\} - 1)] \\
  &=& \left[(b+d+f)^2 - 2(b+d+f) + 1\right] + eb + [af - \min\{a, f, k\}].
\end{eqnarray*}
The first bracket is a perfect square, and the second bracket is
always nonnegative because $a$ and $f$ are nonnegative integers.
Therefore, the bound is indeed nondecreasing in $c$.  By symmetry, the
same is true for $d$, so we can substitute $c = d = 1$ in the bound to
obtain:
\begin{eqnarray*}
  \frac{e(H_k')}{v(H_k')-2} 
  &\leq& \frac{ab + ef + (b+f) + (a+e) + 1 + (\min\{a, f, k\} - 1)}{a+b+1+1+e+f-2} \\
  &=& 1 + \frac{ab + ef + \min\{a, f, k\}}{a+b+e+f} \\
  &\leq& 1 + \frac{ab + ef + \min\big\{\frac{a+f}{2}, k\big\}}{a+b+e+f}.
\end{eqnarray*}
Let $x = (a+b)/2$.  We apply a standard ``smoothing'' technique: note
that if we replace $a = b = x$, the denominator does not change, but
the numerator increases by
\begin{eqnarray*}
  && \left[\left(\frac{a+b}{2}\right)^2 + \min\left\{\frac{\frac{a+b}{2} + f}{2}, k\right\}\right] - \left[ab + \min\left\{\frac{a+f}{2}, k\right\}\right] \\
  &=& \left(\frac{a-b}{2}\right)^2 + \min\left\{\frac{\frac{a+b}{2} + f}{2}, k\right\} - \min\left\{\frac{a+f}{2}, k\right\} \\
  &\geq& \left(\frac{a-b}{2}\right)^2 + \min\left\{\frac{\frac{a+b}{2} + f}{2} - \frac{a+f}{2}, k-k \right\} \\
  &=& \left(\frac{a-b}{2}\right)^2 + \min\left\{\frac{b-a}{4}, 0\right\},
\end{eqnarray*}
which is $\geq 0$ because $a$ and $b$ are integers.  Similarly,
letting $y = (e+f)/2$, we have an analogous result for $e$ and $f$.
Therefore,
\begin{equation}
  \label{ineq:ktt-u-b-e-o:1}
  \frac{e(H_k')}{v(H_k')-2} \leq 1 + \frac{x^2 + y^2 + \min\left\{\frac{x+y}{2}, k\right\}}{2x+2y}.
\end{equation}
The bound is symmetric with respect to $x$ and $y$, so we may assume
without loss of generality that $x \leq y$.  Our next claim is that
increasing $x$ to equal $y$ only increases the bound further.  That
is, we aim to prove the following:
\begin{equation}
  \label{ineq:ktt-u-b-e-o:2}
  \frac{x^2 + y^2 + \min\left\{\frac{x+y}{2}, k\right\}}{2x+2y} \leq \frac{y^2 + y^2 + \min\{y, k\}}{2y+2y}.
\end{equation}
Cross-multiplying to clear denominators, this is equivalent to
\begin{displaymath}
  4 x^2 y + 4y^3 + 4y \min\left\{\frac{x+y}{2}, k\right\} \leq 4xy^2 + 4y^3 + (2x+2y)\min\{y, k\}.
\end{displaymath}
Rearranging terms and using the substitution $\min\{(x+y)/2, k\} \leq
\min\{y, k\}$ (we assumed $x \leq y$), it suffices to show
\begin{equation}
  \label{ineq:ktt-u-b-e-o:3}
  0 \leq (y-x)(4xy - 2 \min\{y, k\}).
\end{equation}
Yet $x = (a+b)/2$ and $u \in U_1 \Rightarrow a \geq 1$, so we must
have $x \geq 1/2$.  This, along with $y \geq x$ and $\min\{y, k\} \leq
y$, establishes \eqref{ineq:ktt-u-b-e-o:3}, and hence
\eqref{ineq:ktt-u-b-e-o:2}.  Applying that to \eqref{ineq:ktt-u-b-e-o:1}, we
obtain
\begin{displaymath}
  \frac{e(H_k')}{v(H_k')-2} \leq 1 + \frac{2y^2 + \min\{y, k\}}{4y}.
\end{displaymath}
If $y \leq k$, then $\min\{y, k\} = y$, so the bound simplifies into
$1 + \frac{y}{2} + \frac{1}{4}$, which is clearly increasing in $y$.
Hence, we may increase $y$ up to $k$ and only increase the bound
further.  Without loss of generality, we may now assume that $y \geq
k$.

Yet if $y \geq k$, then $\min\{y, k\} = k$, so we have the bound
\begin{displaymath}
  \frac{e(H_k')}{v(H_k')-2} \leq 1 + \frac{2y^2 + k}{4y},
\end{displaymath}
and we will show that it is increasing in $y$.  Taking the derivative
of the bound, we obtain
\begin{displaymath}
  \frac{d}{dy} = \frac{1}{2} - \frac{k}{4y^2},
\end{displaymath}
which is $\geq 0$ since we assumed that $y \geq k$, and $k \geq 1$.
The maximum possible value of $y = (e+f)/2$ is $\lfloor t/2 \rfloor$,
so we may increase $y$ to this value and find
\begin{displaymath}
  \frac{e(H_k')}{v(H_k')-2} \leq 1 + \frac{2\left\lfloor \frac{t}{2} \right\rfloor^2 + k}{4 \left\lfloor \frac{t}{2} \right\rfloor}.
\end{displaymath}
The final bound is precisely $e(H_k)/(v(H_k)-2)$, so we are done.  \hfill $\Box$

\vspace{4mm}

\noindent {\bf Proof of Lemma
  \ref{lem:ktt-upper-balanced-extension-odd} for $\boldsymbol{\lfloor t/2 \rfloor
  < k \leq 2 \lfloor t/2 \rfloor}$.}\, The proof in this regime is very
similar to the previous one, but we write it out for completeness.
Consider a consecutive pair $(H_k, H_{k+1})$.  Let $\ell = k - \lfloor
t/2 \rfloor$.  By the construction, $H_k$ has the following structure.
The vertex set is partitioned into $\{V_i\}_i^6$, with $|V_3| =
|V_4| = 1$ and all other $|V_i| = \lfloor t/2 \rfloor$.  The pairs
$(V_1, V_2)$ and $(V_5, V_6)$ are each complete bipartite graphs,
the vertex in $V_3$ is adjacent to all of $V_2 \cup V_4 \cup V_6$, 
the vertex in $V_4$ is adjacent to all of $V_1 \cup V_3 \cup V_5$, 
there is a perfect $\lfloor t/2 \rfloor$-edge matching between $V_1$
and $V_6$, and there is an $(\ell-1)$-edge matching between $V_5$ and
$V_2$.  There are no other edges.  Also, there is a pair of vertices
$u \in V_5$, $v \in V_2$, not involved in the $(\ell-1)$-edge
matching, at which the addition of an edge creates $H_{k+1}$.  Now
consider any family of subsets $U_i \subset V_i$ such that $u \in U_5$
and $v \in U_2$.  Let $H_k'$ be the subgraph of $H_k$ induced by $\cup
U_i$.  We must show that $e(H_k')/(v(H_k')-2) \leq e(H_k)/(v(H_k)-2)$.

For brevity, let $a = |U_1|$, $b = |U_2|$, $c = |U_3|$, $d = |U_4|$,
$e = |U_5|$, and $f = |U_6|$.  Since the edges between $U_5$ and $U_2$
form a matching of at most $\ell-1$ edges which does not involve $u
\in U_5$ or $v \in U_2$, there can be at most $\min\{b-1, e-1,
\ell-1\} = \min\{b, e, \ell\} - 1$ edges there.  On the other hand,
the edges between $U_1$ and $U_6$ also are a matching, and we may
bound their number by $\min\{a, f\}$.  Therefore,
\begin{displaymath}
  \frac{e(H_k')}{v(H_k')-2} \leq \frac{ab + ef + c(b+f) + (a+e)d + cd + \min\{a, f\} + (\min\{b, e, \ell\} - 1)}{a+b+c+d+e+f-2}.
\end{displaymath}
As in the proof of this lemma for $k \leq \lfloor t/2 \rfloor$, our
first step is to show that the bound is nondecreasing in $c$ by
calculating its partial derivative.  The calculations are almost
identical, and we arrive at
\begin{displaymath}
  \text{numerator of } \frac{\partial}{\partial c} = \left[(b+d+f)^2 - 2(b+d+f) + 1\right] + [af - \min\{a, f\}] + [eb - \min\{b, e, \ell\}].
\end{displaymath}
The first bracket is a perfect square, and each of the other brackets
is nonnegative because all variables are nonnegative integers.  Hence
we may set $c=1$ (its maximum value) and increase the bound.  The same
holds for $d$ by symmetry.  Therefore,
\begin{eqnarray*}
  \frac{e(H_k')}{v(H_k')-2} 
  &\leq& \frac{ab + ef + (b+f) + (a+e) + 1 + \min\{a, f\} + (\min\{b, e, \ell\} - 1)}{a+b+1+1+e+f-2} \\
  &=& 1 + \frac{ab + ef + \min\{a, f\} + \min\{b, e, \ell\}}{a+b+e+f} \\
  &\leq& 1 + \frac{ab + ef + \frac{a+f}{2} + \min\{\frac{b+e}{2}, \ell\}}{a+b+e+f}.
\end{eqnarray*}
We again will use a ``smoothing'' technique.  Let $x = (a+b)/2$, and
consider the effect of setting $a = b = x$.  The denominator will
remain invariant, but the numerator will increase by
\begin{eqnarray*}
  && \left[\left(\frac{a+b}{2}\right)^2 + \frac{\frac{a+b}{2} + f}{2} + \min\left\{\frac{\frac{a+b}{2} + e}{2}, \ell\right\}\right] - \left(ab + \frac{a+f}{2} + \min\left\{\frac{b+e}{2}, \ell\right\}\right) \\
  &=& \left(\frac{a-b}{2}\right)^2 + \frac{b-a}{4} + \left[\min\left\{\frac{\frac{a+b}{2} + e}{2}, \ell\right\} - \min\left\{\frac{b+e}{2}, \ell\right\}\right] \\
  &\geq& \left(\frac{a-b}{2}\right)^2 + \frac{b-a}{4} + \min\left\{ \frac{\frac{a+b}{2} + e}{2} - \frac{b+e}{2}, \ell - \ell \right\} \\
  &=& \left(\frac{a-b}{2}\right)^2 + \frac{b-a}{4} + \min\left\{ \frac{a-b}{4}, 0 \right\},
\end{eqnarray*}
which is always $\geq 0$ since $a$ and $b$ are integers.  A similar results holds if we let $y = (e+f)/2$.  Therefore, 
\begin{equation}
  \label{ineq:ktt-u-b-e-o:11}
  \frac{e(H_k')}{v(H_k')-2} \leq 1 + \frac{x^2 + y^2 + \frac{x+y}{2} + \min\left\{\frac{x+y}{2}, \ell\right\}}{2x+2y}.
\end{equation}
The bound is symmetric with respect to $x$ and $y$, so we may assume
without loss of generality that $x \leq y$.  Our next claim is that
increasing $x$ to equal $y$ only increases the bound further.  That
is, we aim to prove the following:
\begin{equation}
  \label{ineq:ktt-u-b-e-o:12}
  \frac{x^2 + y^2 + \frac{x+y}{2} + \min\left\{\frac{x+y}{2}, \ell\right\}}{2x+2y} \leq \frac{y^2 + y^2 + y + \min\{y, \ell\}}{2y+2y}.
\end{equation}
Cross-multiplying to clear denominators, this is equivalent to
\begin{displaymath}
  4 x^2 y + 4y^3 + 2y(x+y) + 4y \min\left\{\frac{x+y}{2}, \ell\right\} \leq 4xy^2 + 4y^3 + (2x+2y)y + (2x+2y)\min\{y, \ell\}.
\end{displaymath}
Rearranging terms and using the substitution $\min\{(x+y)/2, \ell\} \leq \min\{y, \ell\}$ (we assumed $x \leq y$), it suffices to show
\begin{equation}
  \label{ineq:ktt-u-b-e-o:13}
  0 \leq (y-x)(4xy - 2\min\{y, \ell\}).
\end{equation}
Yet $x = (a+b)/2$ and $v \in U_2 \Rightarrow b \geq 1$, so we must
have $x \geq 1/2$.  This, along with $y \geq x$ and $\min\{y, \ell\} \leq
y$, establishes \eqref{ineq:ktt-u-b-e-o:13}, and hence
\eqref{ineq:ktt-u-b-e-o:12}.  Applying that to \eqref{ineq:ktt-u-b-e-o:11}, we
obtain
\begin{displaymath}
  \frac{e(H_k')}{v(H_k')-2} \leq 1 + \frac{2y^2 + y + \min\{y, \ell\}}{4y}.
\end{displaymath}
If $y \leq \ell$, then $\min\{y, \ell\} = y$ and the bound becomes $1
+ \frac{y}{2} + \frac{1}{2}$, which is clearly increasing in $y$.
Hence, we may increase $y$ up to $\ell$ and only increase the bound
further.  Without loss of generality, we may now assume that $y \geq
\ell$.

Yet if $y \geq \ell$, then $\min\{y, \ell\} = \ell$, so we have the bound
\begin{displaymath}
  \frac{e(H_k')}{v(H_k')-2} \leq 1 + \frac{2y^2 + y + \ell}{4y},
\end{displaymath}
and we will show that it is increasing in $y$.  Taking the derivative
of the bound, we obtain
\begin{displaymath}
  \frac{d}{dy} = \frac{1}{2} - \frac{\ell}{4y^2},
\end{displaymath}
which is $\geq 0$ since we assumed that $y \geq \ell$, and $\ell \geq 1$.
The maximum possible value of $y = (e+f)/2$ is $\lfloor t/2 \rfloor$,
so we may increase $y$ to this value and find
\begin{displaymath}
  \frac{e(H_k')}{v(H_k')-2} \leq 1 + \frac{2\left\lfloor \frac{t}{2} \right\rfloor^2 + \left\lfloor \frac{t}{2} \right\rfloor + \ell}{4 \left\lfloor \frac{t}{2} \right\rfloor}.
\end{displaymath}
The final bound is precisely $e(H_k)/(v(H_k)-2)$, so we are done.  \hfill $\Box$

}  

\vspace{4mm}

\noindent {\bf Proof of Lemma
  \ref{lem:ktt-upper-balanced-extension-odd} for $\boldsymbol{k > 2 \lfloor t/2
  \rfloor}$.}\, Consider a consecutive pair $(H_k, H_{k+1})$.  By
construction, $H_k$ has the following structure.  The vertex set is
partitioned into $\{V_i\}_1^6$, with $|V_3| = |V_4| = 1$ and all
other $|V_i| = \lfloor t/2 \rfloor$.  The pairs $(V_1, V_2)$ and
$(V_5, V_6)$ are each complete bipartite graphs, 
the vertex in $V_3$ is adjacent to all of $V_2 \cup V_4 \cup V_6$, 
the vertex in $V_4$ is adjacent to all of $V_1 \cup V_3 \cup V_5$, 
there is a perfect
$\lfloor t/2 \rfloor$-edge matching between $V_1$ and $V_6$, and
another perfect matching between $V_5$ and $V_2$.  There may be some
more edges as well between $V_1$ and $V_6$ or between $V_5$ and $V_2$,
but not all such edges are present: without loss of generality, let us
suppose that there are two vertices $u \in V_1$ and $v \in V_6$ such
that there is no edge between $u$ and $v$.  There are no more edges in
the entire graph.  Also, $H_{k+1}$ is obtained from $H_k$ by 
adding the edge joining $u$ and $v$.  Now, consider any
family of subsets $U_i \subset V_i$ such that $u \in U_1$ and $v \in
U_6$.  Let $H_k'$ be the subgraph of $H_k$ induced by $\cup U_i$.  We
must show that $e(H_k')/(v(H_k')-2) \leq e(H_k)/(v(H_k)-2)$.

For brevity, let $a = |U_1|$, $b = |U_2|$, $c = |U_3|$, $d = |U_4|$,
$e = |U_5|$, and $f = |U_6|$.  Let $E$ be the number of edges in $H_k$
between $U_1$ and $U_6$ or between $U_5$ and $U_2$.  Then
\begin{equation}
  \label{ineq:lemB5-1}
  \frac{e(H_k')}{v(H_k')-2} = \frac{ab + ef + c(b+f) + (a+e)d + cd + E}{a+b+c+d+e+f-2}.
\end{equation}
Next, recall that $H_k$ contained a perfect $\lfloor t/2 \rfloor$-edge
matching between $V_1$ and $V_6$, and between $V_5$ and $V_2$.  The
maximum number of edges of these matchings that are included in $E$
(i.e., go between $U_1$ and $U_6$, or between $U_5$ and $U_2$) is
$\min\{a, f\} + \min\{b, e\} \leq (a+f+b+e)/2$.  Therefore, the number
of edges in $H_k$ between $V_1$ and $V_6$ or between $V_5$ and $V_2$
is at least $E + 2 \big\lfloor \frac{t}{2} \big\rfloor - \frac{a+b+e+f}{2}$.
The rest of the edges in $H_k$ are easy to count: $(V_1, V_2)$ and
$(V_5, V_6)$ are complete bipartite subgraphs $K_{\lfloor t/2 \rfloor,
  \lfloor t/2 \rfloor}$, the vertex in $V_3$ is adjacent to all of
$V_2 \cup V_4 \cup V_6$, and the vertex in $V_4$ is adjacent to all of
$V_1 \cup V_3 \cup V_5$.  Therefore,
\begin{equation}
  \label{ineq:lemB5-2}
  \frac{e(H_k)}{v(H_k)-2} \geq \frac{2 \left\lfloor \frac{t}{2} \right\rfloor^2 + \left[4 \left\lfloor \frac{t}{2} \right\rfloor + 1\right] + \big[E + 2 \left\lfloor \frac{t}{2} \right\rfloor - \frac{a+b+e+f}{2} \big] }{2t-2}.
\end{equation}%
\shortversion{%
  The result follows by proving that the right hand side of
  \eqref{ineq:lemB5-1} is at most the right hand side of
  \eqref{ineq:lemB5-2}.  The full version of this paper contains the
  details.  \hfill $\Box$
}{%
Identifying the denominator as $4 \lfloor t/2 \rfloor$ since $t$ is odd, we may simplify the bound:
\begin{displaymath}
  \frac{e(H_k)}{v(H_k)-2} \geq \frac{1}{2} \left\lfloor \frac{t}{2} \right\rfloor + \frac{3}{2} + \frac{1 + E - \frac{a+b+e+f}{2} }{4 \lfloor t/2 \rfloor},
\end{displaymath}
and so it suffices to show that
\begin{displaymath}
  \frac{ab + ef + c(b+f) + (a+e)d + cd + E}{a+b+c+d+e+f-2} \leq \frac{1}{2} \left\lfloor \frac{t}{2} \right\rfloor + \frac{3}{2} + \frac{1 + E - \frac{a+b+e+f}{2} }{4 \lfloor t/2 \rfloor}.
\end{displaymath}
Observe that $E$ appears in the numerator on both sides.  Since the
denominator on the left hand side is always less than or equal to the
denominator on the right hand side, we may replace $E$ (on both sides) by
its maximum possible value and only sharpen the inequality.  This
maximum is $(|U_1| |U_6|-1) + (|U_5| |U_2|) = af-1+be$, because we
assumed that $H_k$ had no edge between $u \in U_1$ and $v \in U_6$.
Thus it suffices to show that
\begin{displaymath}
  \frac{ab + ef + c(b+f) + (a+e)d + cd + (af-1+be)}{a+b+c+d+e+f-2} \leq \frac{1}{2} \left\lfloor \frac{t}{2} \right\rfloor + \frac{3}{2} + \frac{1 + (af-1+be) - \frac{a+b+e+f}{2} }{4 \lfloor t/2 \rfloor}.
\end{displaymath}
Simplifying the numerator on each side, this is equivalent to
\begin{equation}
  \label{ineq:ktt-u-b-e-o:21}
  \frac{(a+c+e)(b+d+f)-1}{a+b+c+d+e+f-2} \leq \frac{1}{2} \left\lfloor \frac{t}{2} \right\rfloor + \frac{3}{2} + \frac{(a+e)(b+f) - (ab+ef) - \frac{a+b+e+f}{2} }{4 \lfloor t/2 \rfloor}.
\end{equation}
Note that the variables $c$ and $d$ appear only on the left hand side.  We
will show that the left hand side is nondecreasing in each variable, which
will then allow us to replace each of them by their maximum value
of 1 and only sharpen the inequality.  Let us verify this by taking
the partial derivative of the left hand side with respect to $c$:
\begin{eqnarray*}
  \text{numerator of } \frac{\partial}{\partial c} &=& (a+b+c+d+e+f-2)(b+d+f) - [(a+c+e)(b+d+f)-1] \\
  &=& (b+d+f)^2 - 2(b+d+f) + 1.
\end{eqnarray*}
This is a perfect square, so it is $\geq 0$ as claimed.  By symmetry,
a similar result holds for $d$; therefore
\begin{displaymath}
  \frac{(a+c+e)(b+d+f)-1}{a+b+c+d+e+f-2} \leq \frac{(a+1+e)(b+1+f)-1}{a+b+1+1+e+f-2} = 1 + \frac{(a+e)(b+f)}{a+b+e+f},
\end{displaymath}
and applying this to \eqref{ineq:ktt-u-b-e-o:21}, we see that it suffices to show
\begin{equation}
  \label{ineq:ktt-u-b-e-o:22}
  \frac{(a+e)(b+f)}{a+b+e+f} \leq \frac{1}{2} \left\lfloor \frac{t}{2} \right\rfloor + \frac{1}{2} + \frac{(a+e)(b+f) - (ab+ef) - \frac{a+b+e+f}{2} }{4 \lfloor t/2 \rfloor}.
\end{equation}
Our next step is to apply a ``smoothing'' argument to push $(a+e)$ and
$(b+f)$ closer together while preserving their sum.  We claim that
this will only sharpen the inequality.  By symmetry, we may assume
without loss of generality that $a+e \leq b+f$.  Let $\Delta = [(b+f)
- (a+e)]/4$, and consider the effect on \eqref{ineq:ktt-u-b-e-o:22} of
increasing each of $a$ and $e$ by $\Delta$, while decreasing each of
$b$ and $f$ by $\Delta$.  This will bring $(a+e)$ and $(b+f)$ together
to their average, which by convexity, increases their product
$(a+e)(b+f)$.  This product appears in the numerator on both sides,
and since the denominator on the left is at most the denominator on
the right, this simultaneous change can only sharpen the inequality.
The only other term that will change is $ab+ef$, so it remains to
verify that it increases under this transformation:
\begin{displaymath}
  (a+\Delta)(b-\Delta) + (e+\Delta)(f-\Delta) = (ab+ef) + (b+f-a-e)\Delta - 2\Delta^2 = (ab+ef) + 2\Delta^2.
\end{displaymath}
This proves our claim that we may ``smooth'' $(a+e)$ and $(b+f)$
together and only sharpen the inequality.  Therefore, we may now
assume that $a+e = b+f$ for the remainder of the proof; let $x =
(a+e)/2$.  

Rewriting the right hand side of \eqref{ineq:ktt-u-b-e-o:22}, it suffices to show
\begin{displaymath}
  \frac{(a+e)(b+f)}{a+b+e+f} \leq \frac{1}{2} \left\lfloor \frac{t}{2} \right\rfloor + \frac{1}{2} + \frac{af+be - \frac{a+b+e+f}{2} }{4 \lfloor t/2 \rfloor}.
\end{displaymath}
Substituting $a+e = 2x = b+f$ into the inequality, this is equivalent to
\begin{displaymath}
  \frac{(2x)(2x)}{2x+2x} \leq \frac{1}{2} \left\lfloor \frac{t}{2} \right\rfloor + \frac{1}{2} + \frac{af+be - \frac{2x+2x}{2} }{4 \lfloor t/2 \rfloor}.
\end{displaymath}
Simplifying each side of the inequality, this is equivalent to
\begin{equation}
  \label{ineq:ktt-u-b-e-o:23}
  x \leq \frac{1}{2} \left\lfloor \frac{t}{2} \right\rfloor + \frac{1}{2} + \frac{af+be - 2x}{4 \lfloor t/2 \rfloor},
\end{equation}

Observe that we always have $2x \leq 2 \lfloor t/2 \rfloor$ because $0
\leq a, b, e, f \leq \lfloor t/2 \rfloor$; hence $\frac{1}{2} +
\frac{af+be - 2x}{4 \lfloor t/2 \rfloor} \geq 0$.  Therefore, if $x$
itself is already $\leq \frac{1}{2} \big\lfloor \frac{t}{2}
\big\rfloor$, then we are already done.

It remains to consider the case when $2x > \lfloor t/2 \rfloor$, which
we will assume for the remainder of this proof.  Let us replace the
$af+be$ term by a function of $x$, so that the entire inequality will
be in terms of $x$ and $t$ only.  That is, we wish to find the minimum
value of $af+be$ when $a, b, e, f \in [0, \lfloor t/2 \rfloor]$,
subject to the constraints $a+e = 2x = b+f$.  Writing $e = 2x - a$ and
$f = 2x - b$, we aim to minimize
\begin{displaymath}
  af + be = a(2x - b) + b(2x - a) = 2x(a+b) - 2ab.
\end{displaymath}
For fixed $a+b$, the right hand side is minimized when $a=b$, so we
may assume that is the case; let $z = a = b$.  We are then trying to
minimize $2x(2z) - 2z^2$ for fixed $x$.  The variable $z$ is
constrained to the interval $[2x - \lfloor t/2 \rfloor, \lfloor t/2
\rfloor]$, because $0 \leq a \leq \lfloor t/2 \rfloor$ and $0 \leq e =
2x - a \leq \lfloor t/2 \rfloor$ and we assumed $2x > \lfloor t/2
\rfloor$.  Since $2x(2z) - 2z^2$ is quadratic in $z$ with negative
leading coefficient, it is minimized at one of the endpoints of $z$'s
interval.  It turns out that it takes that same value $2 \lfloor t/2
\rfloor (2x - \lfloor t/2 \rfloor)$ at each endpoint, so we conclude
that $af+be \geq 2 \lfloor t/2 \rfloor (2x - \lfloor t/2 \rfloor)$.

Substituting this into \eqref{ineq:ktt-u-b-e-o:23}, it now suffices to
show
\begin{displaymath}
  x \leq \frac{1}{2} \left\lfloor \frac{t}{2} \right\rfloor + \frac{1}{2} + \frac{2 \lfloor t/2 \rfloor (2x - \lfloor t/2 \rfloor) - 2x}{4 \lfloor t/2 \rfloor},
\end{displaymath}
Rearranging terms, this is equivalent to
\begin{displaymath}
  \frac{2x}{4 \lfloor t/2 \rfloor} \leq \frac{1}{2},
\end{displaymath}
which is true because $x = (a+e)/2$ and $a, e \leq \lfloor t/2 \rfloor$, so we are done. \hfill $\Box$

} 

\subsection{Inequalities}

\shortversion{The following inequalities are proved in the full version of this paper.}{
For the reader's convenience, we reproduce the definitions of $s$ and
$\theta$:
\begin{displaymath}
  s = \lfloor \log_r [(r-1)t + 1] \rfloor, \quad\quad\quad \theta = \frac{r^s (2t-2) + 2}{r^s (t^2 - s) + \frac{r^s - 1}{r-1}}.
\end{displaymath}

}  

\begin{inequality}
  \label{ineq:ktt-upper-t/2-codegree}
  Suppose that $t \geq 3$ and $r \geq 2$.  Then $-\theta >
  -\frac{2}{t}$.
\end{inequality}

\shortversion{}{\noindent {\bf Proof.}\, Clearing the denominators, this is equivalent
to
\begin{displaymath}
  r^s (2t^2 - 2t) + 2t < r^s (2t^2 - 2s) + 2 \frac{r^s - 1}{r-1}.
\end{displaymath}
Rearranging terms, this is equivalent to
\begin{displaymath}
  2r^s (s-t) + 2t < 2 \frac{r^s - 1}{r-1}.
\end{displaymath}
Since the right hand side is always positive, it suffices to show that the
left hand side is $\leq 0$.  Note that if $s \leq t/2$, then the left hand side
is indeed $\leq t(-r^s + 2) \leq 0$ since $r \geq 2$ and $s \geq 1$.
This happens whenever $t \geq 4$ by Lemma \ref{lem:s-leq-t/2}, and
whenever $t=3$ and $r \geq 3$ because Lemma \ref{lem:s-dec-in-r} shows
that then $s \leq \lfloor \log_3 (2t+1) \rfloor = 1 \leq t/2$.  The
only remaining case is $t=3$, $r=2$, in which case $s = 2$, and one
may manually verify that the left hand side is precisely $-2 \leq
0$. \hfill $\Box$

\vspace{4mm}

}  

\begin{inequality}
  \label{ineq:ktt-upper-pos-power}
  For any $t \geq 3$ and $r \geq 2$, if $p \gg n^{-\theta}$, then
  $n^{2t-2} p^{t^2 - s - 1}$ is a positive power of $n$.
\end{inequality}

\shortversion{}{\noindent {\bf Proof.}\, The proof of this inequality is nearly
identical to the proof of Inequality \ref{ineq:kt-upper-pos-power}.
Cross-multiplying to clear the denominators, it suffices to show
\begin{equation}
  \label{ineq:i:ktt-u-p-p:1}
  (2t-2)\left[ r^s (t^2 - s) + \frac{r^s - 1}{r-1} \right] > (t^2 - s - 1) \left[r^s (2t-2) + 2\right].
\end{equation}
The right hand side of \eqref{ineq:i:ktt-u-p-p:1} is equal to 
\begin{displaymath}
  (t^2 - s - 1) \left[r^s (2t-2) + 2\right] = (t^2 - s) r^s(2t-2) - r^s (2t-2) + (t^2 - s - 1)2,
\end{displaymath}
so \eqref{ineq:i:ktt-u-p-p:1} is equivalent to
\begin{displaymath}
  (2t-2) \frac{r^s - 1}{r-1} > - r^s (2t-2) + (t^2 - s - 1) 2,
\end{displaymath}
and by rearranging terms and ignoring the $s$ term appearing in the
parentheses on the right hand side, it suffices to establish
\begin{displaymath}
  (2t-2)\frac{r^{s+1} - 1}{r-1} \geq (t^2 - 1)2.
\end{displaymath}
Dividing through by $2(t-1)$, we see that it suffices to show that
\begin{displaymath}
  \frac{r^{s+1} - 1}{r-1} \geq t+1,
\end{displaymath}
which was already shown at the end of the proof of Inequality
\ref{ineq:kt-upper-pos-power}.  \hfill $\Box$

\vspace{4mm}

}  

\begin{inequality}
  \label{ineq:ktt-upper-neg-power}
  For any $t \geq 3$ and $r \geq 2$, if $p \ll n^{-\theta}$, then
  $n^{2t-2} p^{t^2 - s}$ is a negative power of $n$.
\end{inequality}

\shortversion{}{\noindent {\bf Proof.}\, The proof of this inequality is nearly
identical to the proof of Inequality \ref{ineq:kt-upper-neg-power}.
We must establish the following inequality:
\begin{displaymath}
  (2t-2)\left[ r^s (t^2 - s) + \frac{r^s - 1}{r-1} \right] < (t^2 - s) \left[r^s (2t-2) + 2\right].
\end{displaymath}
This is equivalent to
\begin{displaymath}
  (2t-2)\frac{r^s - 1}{r-1} < (t^2 - s) 2.
\end{displaymath}
By definition of $s$, $r^s \leq (r-1)t+1$.  Therefore, the left hand
side is at most $(2t-2)t$, and it suffices to show that
\begin{displaymath}
  2t^2 - 2t < 2t^2 - 2s,
\end{displaymath}
i.e., to show that $s < t$.  But Lemma \ref{lem:s-dec-in-r}
established that if $r \geq 2$, then $s \leq \lfloor \log_2 (t+1)
\rfloor$, and this is indeed less than $t$ for all $t \geq 3$.  \hfill
$\Box$

}  

\end{document}